\documentclass{siamart250106}

\usepackage{amssymb,amsmath}
\usepackage{bbm}
\usepackage{graphicx}
\usepackage{caption}
\usepackage{subcaption}
\usepackage{hyperref}
\usepackage{booktabs}
\usepackage{multirow}
\usepackage{enumerate}

\hypersetup{
     colorlinks   = true,
     citecolor    = blue
}

\usepackage{mathtools}
\DeclarePairedDelimiter{\ceil}{\lceil}{\rceil}

\newcommand{\R}{\mathbb R}

\newcommand{\N}{\mathbb N}

\DeclareMathOperator{\sgn}{sgn}
\DeclareMathOperator*{\argmax}{arg\,max}
\DeclareMathOperator*{\argmin}{arg\,min}

\DeclareFontFamily{U}{mathx}{}
\DeclareFontShape{U}{mathx}{m}{n}{<-> mathx10}{}
\DeclareSymbolFont{mathx}{U}{mathx}{m}{n}
\DeclareMathAccent{\widehat}{0}{mathx}{"70}
\DeclareMathAccent{\widecheck}{0}{mathx}{"71}

\newtheorem{remark}[theorem]{Remark}

\title{An open waveguide with a thin high contrast core layer: asymptotic analysis and inverse detection problem
}

\author{
E. Bonnetier\thanks{Institut Fourier, Universit\'e Grenoble-Alpes, France
    (\email{eric.bonnetier@univ-grenoble-alpes.fr}).}
\and M. Courdurier\thanks{Facultad de Matem\'aticas, Pontificia Universidad Cat\'olica de Chile, Chile (\email{mcourdurier@mat.puc.cl}).}
\and A. Osses\thanks{Departamento de Ingenier\'{\i}a Matem\'atica and Center for Mathematical Modeling, Universidad de Chile, Chile (\email{axosses@dim.uchile.cl}).}
\and F. Triki\thanks{Laboratoire Jean Kuntzmann, Universit\'e Grenoble-Alpes, France (\email{faouzi.triki@univ-grenoble-alpes.fr})}
}

\headers{Open waveguide with thin high contrast core layer}{E. Bonnetier, M. Courdurier, A. Osses, F. Triki}

\begin{document}
\maketitle
\begin{abstract}
We investigate the Helmholtz equation in a two dimensional open waveguide with a thin and high contrast core layer. 
We develop an asymptotic analysis of the Green function of the problem, and through it we identify and characterize the appearance of resonant frequencies. For waves originating outside
of the core, the waveguide response at these resonant frequencies is vastly different than the response at non-resonant frequencies. Using this phenomenon and multifrequency measurements containing the first resonance, we propose, theoretically analyze, and numerically validate a reconstruction algorithm to identify the location, thickness and index of refraction of the core layer.
\end{abstract}

\section{Introduction and Motivation}
\label{sec:introduction}

We consider the solutions of the Helmholtz equation in an open wave guide, that consists
in a thin, unbounbed core layer with a high optical index embedded in a homogeneous background
medium.
We are interested in the identification of the descriptive parameters of the core layer from its response to a localized excitation. 
This is related to some inverse problems appearing in 
seismology for layered media 
and in optical or sound probing of laminated media. 

The method we propose for this identification relies on a remarkable resonance phenomenon~:
When the medium is excited by an external point source over a wide range of frequencies,
the reflection and transmission properties of the core layer
abruptly changes at particular frequencies.
This behavior has been observed experimentally and has been reported in the literature 
under the name of {\em absorption power}~\cite{absorption-power} or 
{\em detuned frequencies} ~\cite{detuned}.
It has also been reported in theoretical studies, see for instance~\cite{Ammari-Triki-2004} 
and~\cite{Bonnetier-Triki-2009}.
An explanation put forward for this change in behavior, is that the field inside the 
core layer become resonant (in the sense of the resonant modes of closed waveguides, this will be made clearer in Remark~\ref{rmk:interpretation})
and are thus able to couple the fields on both sides of the core, inspite of the 
large value of the core index that should make it act as an effective Dirichlet
boundary condition.
This phenomenon is reminiscent of the effect of plasmonic waves in periodic gratings (although a somewhat dual mechanism is at work here), which may greatly enhance the fields near the surface of the grating at certain frequencies~\cite{Quemerais_et_al, Ebbesen}.
Our set up also relates to the work on photonic crystals~\cite{Figotin_Kuchment1,Figotin_Kuchment2,Kuchment2001},
where one conjugates thin structures and high contrast to evidence the possibility of
opening band gaps.

Throughout this paper we will refer to the frequencies at which the device changes its
reflection or transmission properties as resonances, a terminology that will be
justified by our asymptotic analysis.
In order to better understand this phenomenon, we explicitly analyze the Green function
for the Helmholtz equation in an open waveguide. 
Finding such Green function in a medium that contains a non-compact inclusion is not
easy in general since finding the right radiation condition can be a delicate matter (see for instance  \cite{Ciraolo-2009, Ciraolo-Magnanini-2009, Duran, Perthame}).
In the simple geometry we consider, where the index of refraction depends only on the 
transversal variable, we rely on the work of R. Magnanini and F. Santosa 
who derived an explicit expression of the Green function 
by passing to the limit in the height of a slab containing the 
core layer~\cite{Magnanini-Santosa-2000}.  
Uniqueness of the Green function is proved in~\cite{Ciraolo-Magnanini-2009}, 
and a generalized form is given in~\cite{Ciraolo-2009} when the index of refraction 
may be perturbed in the longitudinal direction. 
We also refer to~\cite{Jerez-Nedelec-2012} and the references therein,
who studies open three dimensional waveguides, where the proper radiation conditions are obtained via a limiting absorption principle. 

The current article has two main objectives: 1) We derive an asymptotic form of the Green function
in the case of a thin core layer with a high index of refraction, which explains the appearance
of resonances; ans 2) We take advantage of this resonance phenomenon and of the associated asymptotic
form of the Green function to address the inverse problem of detecting the location of the core layer,
its thickness, and index of refraction from measurements on a bounded receiving surface.

Our work is organized as follows: in Section~\ref{sec:setup}, we introduce 
the Green function~$G$ of the open waveguide we consider, a straight core layer 
of thickness $2h$, parallel to the $z$-axis, in the $(x,z)$ plane, with optical index 
$n_h = \bar{n}/h$, surrounded by a homogeneous medium with fixed optical index. 
We also recall the form of the Green function of the whole plane with constant index, which is represented in terms of the Hankel function.
Our asymptotic results are presented in Section~\ref{sec:asymptotic}~: Theorem~\ref{thm:asym_guided}
describes the limiting
behavior of the part of the Green function~$G$ that corresponds to guided modes, while in
Theorems~\ref{thm:asym_sym} and~\ref{thm:asym_anti}, we give the limiting behavior of its continuous symmetric
and antisymmetric parts, which depend on the fact that the frequency $k$ is resonant or not.
We summarize the asymptotic behavior of the Green function in Theorem~\ref{thm:total_G}~: 
as $h \to 0$, depending on $k$, $G$ tends to the Green function of a half-plane
with constant index and Dirichlet boundary condition, or to the Green function of the whole plane. In the non-resonant case, in the limit, the core behaves as an absorbant boundary, whereas in the resonant case, the core becomes invisible for an observer located on the same side as the source.

In Section~\ref{sec:inverse_problem}, we propose and analyze an algorithm for the inverse problem of recovering the position,  thickness, and optical index of the core layer in the regime $h \to 0$. To take advantage of the resonance phenomenon, we use multifrequency measurements.
In Section \ref{sec:numerical-results}, we report numerical results. They show that our proposed algorithm is capable of successfully recovering the characterizing parameters of the core layer, even in the presence of relatively large noise.
Finally, Section~\ref{sec:asymptotic-formulas} is devoted to proving the asymptotic results of Section~\ref{sec:asymptotic}.

\section{Setup and the Green Functions}
\label{sec:setup}

We consider $\R^2=\{(x,z):x,z\in\R\}$ as an open waveguide that
consists in a thin core layer $\{(x,z)\in\R^2: |x|<h\}$ of width $2h$ embedded
in a homogeneous cladding medium. In this setting we consider the following
Helmholtz equation
\begin{align}\label{eq:Helmholtz1}
\Delta u+k^2n^2(x)u=f, \text{ in } \R^2
\end{align}
where the refraction index depends only in the $x$ coordinate as
\begin{align}\label{eq:Helmholtz2}
n(x,z)=n(x)=
\begin{cases}
n_{cl}&\textnormal{ for } |x|>h\\
n_h=\overline{n}/h&\textnormal{ for } |x|<h,
\end{cases}
\end{align}
and $n_h>n_{cl}$ (the core and cladding indexes of refraction respectively). We assume that
$h$ is small and that $n_h=\overline{n}/h$ is large. Throughout the paper, we assume that $f$ is a Dirac mass located at a point $(x_0,z_0)$ above the core.

\subsection{The Green function $G$}

The Green function for this Helmholtz
equation, with the adequate radiation condition, is studied in
\cite{Magnanini-Santosa-2000}, and it takes the form described below.
The formula is given for $(x,z),(x_0,z_0)\in\R^2$, but for our purpose, we are only
interested in the case when $x_0\geq c$, for some constant $c>0$.

With the convention $i=\sqrt{-1}$ we introduce the following notation. Let
$$d^2=k^2(n_h^2-n_{cl}^2),\quad \lambda=k^2(n_h^2-\beta^2)m \quad k\beta=\sqrt{k^2n_h^2-\lambda},\quad \text{and} \quad Q=\sqrt{\lambda-d^2}.$$

Let $v_s$ and $v_a$, that we call symmetric and antisymmetric modes, be defined by
\begin{align*}
&v_s(x,\lambda)\!=\!\left\{
\begin{array}{cl}
\cos(h\sqrt{\lambda})\cos(Q(|x|-h))
-\frac{\sqrt{\lambda}}{Q}\sin(h\sqrt{\lambda})\sin(Q(|x|-h))
&, |x|>h\\
\cos(x\sqrt{\lambda})&, |x|\le h
\end{array}\right.\\
&v_a(x,\lambda)\!=\!\left\{
\begin{array}{cl}
\!\!\!\sgn(x)\!\!\left[\sin(h\sqrt{\lambda})\cos(Q(|x|-h))\!+\!\frac{\sqrt{\lambda}}{Q}\cos(h\sqrt{\lambda})\sin(Q(|x|-h))\right]
&\!\!\!\!\!\!, |x|>h\\
\sin(x\sqrt{\lambda})&\!\!\!\!\!\!, |x|\le h\\
\end{array}
\right.
\end{align*}

Let $\{\lambda_{s,j}\}_{j=1}^{J_s},\{\lambda_{a,j}\}_{j=1}^{ J_a}$ be the finite set of real roots of the equations:
\begin{align*}
\textnormal{for $\lambda_s$}:\sqrt{d^2-\lambda}-\sqrt{\lambda}\tan(h\sqrt{\lambda})=0,\qquad 
\textnormal{for $\lambda_a$}:\sqrt{d^2-\lambda}+\sqrt{\lambda}\cot(h\sqrt{\lambda})=0,
\end{align*}
which are strictly bounded between 0 and $d^2$ and such that 
\begin{align*}
&v_s(x,\lambda_s)=
\begin{cases}
\cos(h\sqrt{\lambda_s})e^{-\sqrt{d^2-\lambda_s}(|x|-h)}&\textnormal{ if } |x|>h\\
\cos(x\sqrt{\lambda_s})&\textnormal{ if } |x|\le h
\end{cases}\\
&v_a(x,\lambda_a)=
\begin{cases}
\sgn(x)\sin(h\sqrt{\lambda_s})e^{-\sqrt{d^2-\lambda_s}(|x|-h)}&\textnormal{ if } |x|>h\\
\sin(x\sqrt{\lambda_s})&\textnormal{ if } |x|\le h
\end{cases}
\end{align*}

Let $d\rho^s,d\rho^a$ denote the spectral measures:
\begin{align*}
<d\rho^s,\eta>=&\sum_{j}\frac{\sqrt{d^2-\lambda_{s,j}}}{1+h\sqrt{d^2-\lambda_{s,j}}}\eta(\lambda_{s,j})
+ \frac{1}{2\pi}\int_{d^2}^\infty \frac{\sqrt{\lambda-d^2}}{(\lambda-d^2)+d^2\sin^2(h\sqrt{\lambda})}\eta(\lambda) d\lambda
\end{align*}
\begin{align*}
<d\rho^a,\eta>=&\sum_{j}\frac{\sqrt{d^2-\lambda_{a,j}}}{1+h\sqrt{d^2-\lambda_{a,j}}}\eta(\lambda_{a,j})
+ \frac{1}{2\pi}\int_{d^2}^\infty \frac{\sqrt{\lambda-d^2}}{(\lambda-d^2)+d^2\cos^2(h\sqrt{\lambda})}\eta(\lambda) d\lambda.
\end{align*}

Then, the Green function with a point source at $(x_0,z_0)$ is given by
\begin{align}\label{eq:def_G}
G(x,z;x_0,z_0)=\sum_{m\in\{s,a\}}\int v_m(x,\lambda)v_m(x_0,\lambda)\frac{e^{ik\beta|z-z_0|}}{2ik\beta}d\rho^m(\lambda).
\end{align}
More explicitly,
\begin{align*}
G(x,z;x_0,z_0)= G_s(x,z;x_0,z_0)+ G_a(x,z;x_0,z_0),
\end{align*}
where $G_m=G_{m,g}+G_{m,c}$, for $m\in\{s,a\}$, and
\begin{eqnarray*}
{G_{m,g}(x,z;x_0,z_0)}&=&\sum_j v_m(x,\lambda_{m,j})v_m(x_0,\lambda_{m,j})
\frac{e^{ik\beta_{m,j}|z-z_0|}}{2ik\beta_{m,j}}\frac{\sqrt{d^2-\lambda_{m,j}}}{1+h\sqrt{d^2-\lambda_{m,j}}}\\
{G_{s,c}(x,z;x_0,z_0)}&=&\frac{1}{2\pi}\int_{d^2}^\infty v_s(x,\lambda)v_s(x_0,\lambda)
\frac{e^{ik\beta|z-z_0|}}{2ik\beta}\frac{\sqrt{\lambda-d^2}}{(\lambda-d^2)+d^2 \sin^2(h\sqrt{\lambda})} d\lambda\\
{G_{a,c}(x,z;x_0,z_0)}&=&\frac{1}{2\pi}\int_{d^2}^\infty v_a(x,\lambda)v_a(x_0,\lambda)
\frac{e^{ik\beta|z-z_0|}}{2ik\beta}\frac{\sqrt{\lambda-d^2}}{(\lambda-d^2)+d^2 \cos^2(h\sqrt{\lambda})} d\lambda.\\
\end{eqnarray*}
where the subscripts $g$ and $c$ stand for continuous and guided respectively (and $s,a$ stand for symmetric and antisymmetric respectively).
Numerical computations of this Green's function are presented in Section \ref{sec:numerical-results}.

\subsection{The Green function $H$ for the Helmholtz equation in the homogeneous plane}

To interpret the limiting behavior of $G$, the following expression of the Green function for the homogeneous Helmholtz equation in the whole plane is useful.

\begin{lemma}\label{lemma:H} Let $(x_0,z_0)\in\R^2$. Let $H(x,z;x_0,z_0)$ solve the Helmholtz equation in $\R^2$ with an homogeneous index of refraction $n_{cl}$, at frequency $k$, and with a source at $(x_0,z_0)$,
\begin{align*}
    \Delta H(x,z;x_0,z_0)+k^2n_{cl}^2H(x,z;x_0,z_0)=\delta (x-x_0)\delta(z-z_0), \quad x,z\in\R,
\end{align*}       
    considering outgoing Sommerfeld radiation condition. Then
\begin{align*}
H(x,z;x_0,z_0)=&\frac{1}{2\pi}\int_{0}^\infty \cos(k\tau(x-x_0)) 
\frac{e^{ik\sqrt{n^2_{cl}-\tau^2}|z-z_0|}}{i\sqrt{n^2_{cl}-\tau^2}} d\tau.
\end{align*}
\end{lemma}

\begin{proof}
Taking Fourier transform in $x\mapsto \omega$ variable, and denoting by $\hat H$ the Fourier transform of $H$ in the $x$ variable, the equation becomes
\begin{align*}
    \partial_z^2 \hat H +(k^2n^2-\omega^2)\hat H = e^{-i\omega x_0}\delta(z-z_0),
\end{align*}
which then implies that for each $\omega$, and considering the Sommerfeld radiation condition,  
\begin{align*}
    \hat H = e^{-i\omega x_0}\frac{ e^{i\sqrt{k^2n^2-\omega^2}|z-z_0|}}{2i\sqrt{k^2n^2-\omega^2}}.
\end{align*}
Taking inverse Fourier transform in $\omega$ we get
\begin{align*}
    H = \frac{1}{2\pi}\!\!\int_{-\infty}^{\infty}\!\!\!\! e^{i\omega z}e^{-i\omega x_0}\frac{ e^{i\sqrt{k^2n^2-\omega^2}|z-z_0|}}{2i\sqrt{k^2n^2-\omega^2}} d\omega
    = \frac{1}{2\pi}\!\!\int_0^{\infty} \!\!\!\!\cos(\omega(x-x_0))\frac{ e^{i\sqrt{k^2n^2-\omega^2}|z-z_0|}}{i\sqrt{k^2n^2-\omega^2}} d\omega,
\end{align*}
and considering the change of variable $\omega=k\tau$ we conclude
\begin{align*}
    H(x,z;x_0,z_0;k) =& \frac{1}{2\pi}\int_0^{\infty} \cos(k\tau(x-x_0))\frac{ e^{ik\sqrt{n^2-\tau^2}|z-z_0|}}{i\sqrt{n^2-\tau^2}} d\tau.
\end{align*}
\end{proof}

\begin{remark}
    The formula for $H$ is consistent with the expression of $G$ in~\eqref{eq:def_G} when choosing $n_h=n_{cl}$. And as observed in \cite{Ciraolo-Magnanini-2009}, this Green function
satisfies the weaker version of the Sommerfeld radiation condition established by Rellich \cite{Rellich-1943}. In other words, the above expression of $H$ is an alternative to the usual expression given in terms of the Hankel function.
\end{remark}

As in the case of $G$, we can split $H$ in its symmetric and antisymmetric parts with respect to $x$, namely $H=H_s+H_a$ with
\begin{align}
\begin{split}\label{eq:Ha}
    H_s(x,z;x_0,z_0)=&(H(x,z;x_0,z_0)+H(-x,z;x_0,z_0))/2,\\
    H_a(x,z;x_0,z_0)=&(H(x,z;x_0,z_0)-H(-x,z;x_0,z_0))/2,\\
\end{split}
\end{align}
and where the reflection done in the $x$ variable could have been done in the $x_0$ variable instead. Note that $2H_a$ (respectively $2H_s$) is the Green function of the Helmholtz equation in the upper half-plane with Dirichlet (respectively Neumann) boundary condition on $\{x=0\}$.

\section{Asymptotic behavior of the Green function $G$}
\label{sec:asymptotic}

In this section, we study the asymptotic behavior of the various components (guided modes, symmetric continuous, and antisymmetric continuous parts) of $G$. The proofs of the following theorems will be given in Section~\ref{sec:asymptotic-formulas}. We start with guided modes.

\subsection{Asymptotics of the guided component}

\begin{theorem}\label{thm:asym_guided}
Assume that $x_0\geq c>0$.
Then, there exists constants $\kappa_1,\kappa_2>0$ depending only on $k$ and $\overline{n}$, such that for all small enough $h>0$ (in particular requires $c>2h$),
\begin{align}\label{eq:estim_guided}
\|(G_{s,g}+G_{a,g})(\cdot,\cdot;x_0,z_0)\|_{L^\infty(\R^2)}\leq \kappa_1 e^{-\kappa_2\frac{c}{h}}.
\end{align}
\end{theorem}

This result means that, asymptotically, the guided modes vanish very rapidly as $h\to0$ when the source is located outside the core. We proceed to study $G_{s,c}$, the symmetric continuous part of the Green function.

\subsection{Asymptotics of the symmetric continuous component}

The limiting behavior of $G_{s,c}$ depends on whether $\sin(k\overline{n})$ vanishes or not.

\begin{theorem}\label{thm:asym_sym} Assume that $x_0\geq c>0$.
\begin{enumerate}[(i)]
\item\label{enum:sym1} Assume that $\sin(k\overline n)\neq0$. 
\begin{enumerate}
    \item Let $(x,z)\in K$, a compact set contained in $\{x>0, |z-z_0|>0\}$, i.e. $(x,z)$ is located on the same side of the core as the source. Then for $h>0$ small enough
\begin{align}\label{eq:Gsc_nr_+}
\| (G_{s,c}-[H_a+h\Phi_s])(\cdot,\cdot;x_0,z_0)\|_{L^\infty(K)}\leq \kappa h^2.
\end{align}

    \item Let $(x,z)\in K$, a compact set contained in $\{x<0, |z-z_0|>0\}$, i.e. $(x,z)$ is located on the opposite side of the core from the source. Then for $h>0$ small enough
\begin{align}\label{eq:Gsc_nr_-}
\| (G_{s,c}-[-H_a+h\Phi_s])(\cdot,\cdot;x_0,z_0)\|_{L^\infty(K)}\leq \kappa h^2,
\end{align}

\end{enumerate}

\item\label{enum:sym2} Assume that $\sin(k\overline n)=0$. Let $(x,z)\in K$, a compact set contained in $\{|x|>0, |z-z_0|>0\}$. Then for $h>0$ small enough
\begin{align}\label{eq:Gsc_r}
&\| (G_{s,c}-[H_s+h\Psi_s])(\cdot,\cdot;x_0,z_0)\|_{L^\infty(K)}\leq \kappa h^{3/2},
\end{align}

\end{enumerate}
In all these estimates, $\kappa>0$ depends only on $K, k,\overline{n},c$ and $n_{cl}$,
 $H_s$ and $H_a$ are given by \eqref{eq:Ha}, and $\Phi_s$, $\Psi_s$ are defined by
 \begin{eqnarray*}
  {\Phi_s(x,z;x_0,z_0)}
  &=& -\frac{k}{2\pi}\int_{0}^\infty \tau\left(\frac{\cot(k\overline n)}{k\overline n}+1\right)
  \sin(k\tau(|x|+|x_0|))
  \frac{e^{ik\sqrt{n^2_{cl}-\tau^2}|z-z_0|}}{i\sqrt{n^2_{cl}-\tau^2}} d\tau,\\
{\Psi_s(x,z;x_0,z_0)}&=&\frac{k}{2\pi}\int_0^\infty \frac{\tau^2+n_{cl}^2}{2\tau}\sin(k\tau(|x|+|x_0|))
\frac{e^{ik\sqrt{n^2_{cl}-\tau^2}|z-z_0|}}{i\sqrt{n^2_{cl}-\tau^2}} d\tau
- kn_{cl} \frac{e^{ik n_{cl}|z-z_0|}}{8i}.
\end{eqnarray*}

\end{theorem}

\subsection{Asymptotics of the antisymmetric continuous component}

We now turn to the study of $G_{a,c}$, the antisymmetric continuous component of $G$.
In this case, the fact that $\cos(k\overline{n})$ vanishes or not, plays a role.

\begin{theorem}\label{thm:asym_anti} Assume that $x_0\geq c>0$.
\begin{enumerate}[(i)]
\item\label{enum:anti1} Assume that $\cos(k\overline n)\neq0$. Let $(x,z)\in K$, a compact set contained in $\{|x|>0, |z-z_0|>0\}$. Then for $h>0$ small enough
\begin{align}\label{eq:Gsc_r}
&\| (G_{a,c}-[H_a+h\Phi_a])(\cdot,\cdot;x_0,z_0)\|_{L^\infty(K)}\leq \kappa h^2,
\end{align}
\item\label{enum:anti2} Assume that $\cos(k\overline n)=0$. 
\begin{enumerate}
    \item Let $(x,z)\in K$, a compact set contained in $\{x>0, |z-z_0|>0\}$, i.e. $(x,z)$ is located on the same side of the core as the source. Then for $h>0$ small enough
\begin{align}\label{eq:Gsc_nr_+}
\| (G_{a,c}-[H_s+h\Psi_a])(\cdot,\cdot;x_0,z_0)\|_{L^\infty(K)}\leq \kappa h^{3/2}.
\end{align}
    \item Let $(x,z)\in K$, a compact set contained in $\{x<0, |z-z_0|>0\}$, i.e. $(x,z)$ is located on the opposite side of the core from the source. Then for $h>0$ small enough
\begin{align}\label{eq:Gsc_nr_-}
\| (G_{a,c}-[-H_s+h\Psi_a])(\cdot,\cdot;x_0,z_0)\|_{L^\infty(K)}\leq \kappa h^{3/2},
\end{align}
\end{enumerate}
\end{enumerate}
In all these estimates, $\kappa>0$ depends only on $K, k,\overline{n},c$ and $n_{cl}$,
 $H_s$ and $H_a$ are given by \eqref{eq:Ha}, and $\Phi_a$, $\Psi_a$ are defined by
 \begin{eqnarray*}
  {\Phi_a(x,z;x_0,z_0)} 
  =& \frac{k\sgn(xx_0)}{2\pi}\!\!\int_{0}^\infty \!\!\!\tau\left(\frac{\tan(k\overline n)}{k\overline n}\!-\!1\!\right)
  \!\sin(k\tau(|x|+|x_0|))
  \frac{e^{ik\sqrt{n^2_{cl}-\tau^2}|z-z_0|}}{i\sqrt{n^2_{cl}-\tau^2}} d\tau,\\
{\Psi_a(x,z;x_0,z_0)}=&\frac{k\sgn(xx_0)}{2\pi}\int_0^\infty \frac{3\tau^2-n_{cl}^2}{2\tau}\sin(k\tau(|x|+|x_0|))
\frac{e^{ik\sqrt{n^2_{cl}-\tau^2}|z-z_0|}}{i\sqrt{n^2_{cl}-\tau^2}} d\tau\\
&- k\sgn(xx_0)n_{cl} \frac{e^{ik n_{cl}|z-z_0|}}{8i}.
\end{eqnarray*}
\end{theorem}

\subsection{Asymptotic behavior of the complete Green function}

We synthesize the results of the previous theorems in the Theorem below, which sheds light on the remarkable behavior of the core layer as $h$ tends to zero, depending on the frequency.

\begin{theorem}\label{thm:total_G}
Assume that $c > 0$. As $h \to 0$, the behavior of $G_s$ and $G_a$, the symmetric and antisymmetric components of the Green function, and the behaviour of the Green function $G$,  up to order one in $h$, in terms of frequency, is given in the following table.

\setlength\tabcolsep{0pt}
\begin{tabular}{|c|c|c|}
\hline
Frequency & $x > c$ &  $x < -c$ \\
\hline
$\begin{array}{c}
\textrm{non resonant case~:} \\[4pt]
\sin(k \overline{n}) \cos(k \overline{n}) \neq 0
\end{array}$
&
$\begin{array}{c}
G_{s,c}  \sim  H_a + h \Phi_s \\
G_{a,c}  \sim  H_a + h \Phi_a\\
~\\
G \sim 2 H_a+h(\Phi_s+\Phi_a)
\end{array}$
&
$\begin{array}{c}
G_{s,c}  \sim  -H_a + h \Phi_s \\
G_{a,c}  \sim   H_a + h \Phi_a\\
~\\
G \sim 0+ h(\Phi_s+\Phi_a)
\end{array}$
\\
\hline
$\begin{array}{c}
\textrm{sym. resonnant~:} \\[4pt]
\sin(k \overline{n}) = 0
\end{array}$
&
$\begin{array}{c}
G_{s,c}  \sim  H_s + h \Psi_s \\
G_{a,c}  \sim  H_a + h \Phi_a\\
~\\
G \sim H+ h(\Psi_s+\Phi_a)
\end{array}$
&
$\begin{array}{c}
G_{s,c}  \sim  H_s + h \Psi_s \\
G_{a,c}  \sim   H_a + h \Phi_a\\
~\\
G \sim H+ h(\Psi_s+\Phi_a)
\end{array}$
\\
\hline
$\begin{array}{c}
\textrm{Antisym. resonnant~:} \\[4pt]
\cos(k \overline{n}) = 0
\end{array}$
&
$\begin{array}{c}
G_{s,c}  \sim  H_a + h \Phi_s \\
G_{a,c}  \sim  H_s + h \Psi_a\\
~\\
G \sim H+ h(\Phi_s+\Psi_a)
\end{array}$
&
$\begin{array}{c}
G_{s,c}  \sim  -H_a + h \Phi_s \\
G_{a,c}  \sim   -H_s + h \Psi_a\\
~\\
G \sim -H+ h(\Phi_s+\Psi_a)
\end{array}$
\\
\hline
\end{tabular}\\

In particular, at order zero in $h$, the Green function $G$ satisfies,
\begin{center}
\begin{tabular}{|c|c|c|}
\hline
Frequency & ~~~~$x > c$ ~~~~&  ~~~~$x < -c$~~~~ \\
\hline
$\begin{array}{c}
\textrm{Non resonant case~:} \\
\sin(k \overline{n}) \cos(k \overline{n}) \neq 0
\end{array}$
&
$G \sim 2 H_a$
&
$G \sim 0$
\\
\hline
$\begin{array}{c}
\textrm{Symmetric resonant~:} \\
\sin(k \overline{n}) = 0
\end{array}$
&
$G \sim H$
&
$G \sim H$
\\
\hline
$\begin{array}{c}
\textrm{Antisymmetric resonant~:} \\
\cos(k \overline{n}) = 0
\end{array}$
&
$G \sim H$
&
$G \sim -H$ \\
\hline
\end{tabular}
\end{center} 
\end{theorem}

\begin{center}
\begin{figure}
\includegraphics[trim=0 0 0 0, clip, width=13cm]{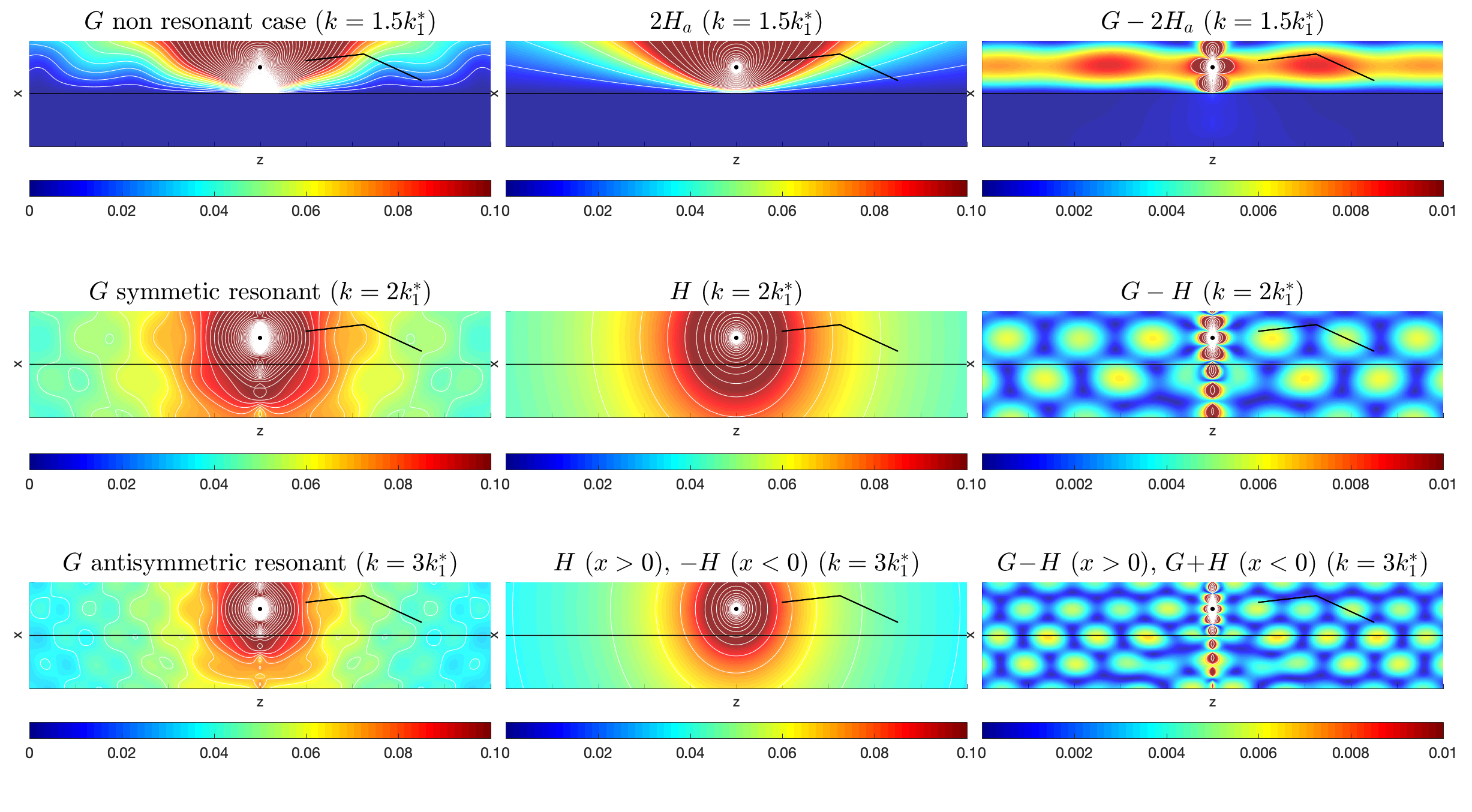}
\caption{Non-resonant (first row) and resonant cases (second and third rows) for an open waveguide with a thin
high contrast core layer. The waveguide parameters are $h=0.005$, $n_{cl}=1$, $n_h=\frac{\pi}{2h}$, $k_1^*=1$. We compare the modulus of the Green functions \eqref{eq:def_G}
(first column), their zero order approximations in $h$,  described in Theorem~\ref{thm:total_G} and written in terms of the Hankel function (second column), and the corresponding error (third column).}
\label{fig:nonresonant_resonant_cases}
\end{figure}
\end{center}

\begin{center}
\begin{figure}
\centering
\includegraphics[trim=0 0 0 0, clip, width=13cm,height=7cm]{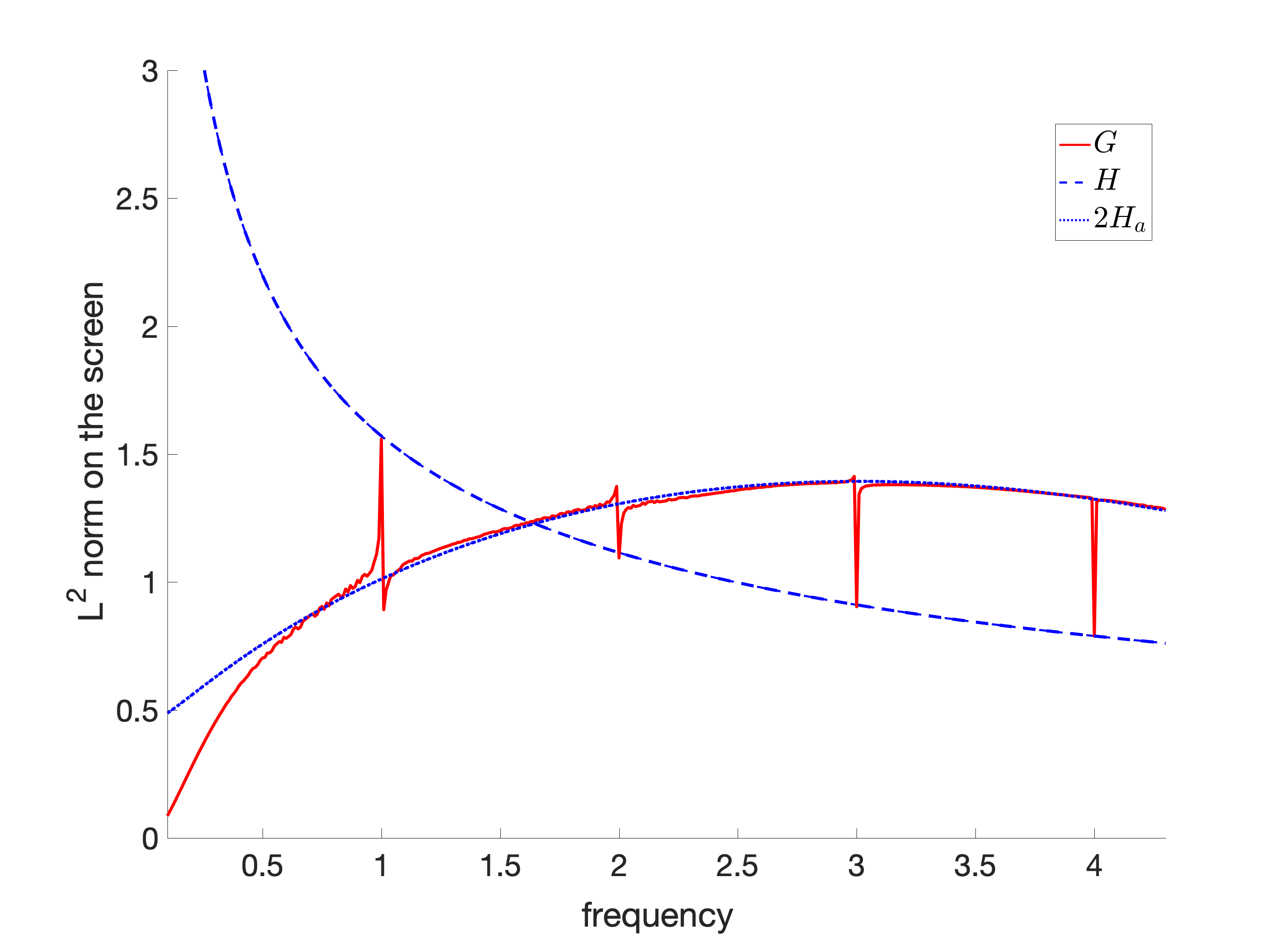}
\caption{Measurements on a screen $S$ for the geometry setting of 
Figure~\ref{fig:geometry-main} for an open waveguide with 
parameters $h=0.005$, $n_{cl}=1$, $n_h=\frac{\pi}{2h}$, $k_1^*=1$. 
Comparison of $||G||_{L^2(S)}$ (red line), $||H||_{L^2(S)}$ (blue dotted line), and $||2H_a||_{L^2(S)}$ (blue line) in terms of frequency. Notice the change of behavior 
(narrow peaks) at the resonant frequencies $k_p^*$, $p=1,2,3,4,\ldots$.}
\label{fig:plotG}
\end{figure}
\end{center}
\begin{remark}\label{rmk:interpretation}
\begin{enumerate}[(i)]
    \item If $\sin(k\overline{n})\cos(k\overline{n})\neq0$. For an observer located above the core (at a point $(x,z)$ with $x>c$) the field is asymptotically the same as if the core layer were acting as a Dirichlet boundary condition, corresponding to an absorbing layer. For an observer located on the opposite side of the core (at a point $(x,z)$ with $x<-c$) the field would simply be zero.
    \item If $\sin(k\overline{n})\cos(k\overline{n})=0$. For an observer located above the core, the field is asymptotically as if the core layer were absent. The same holds for an observer located on the opposite side of the core if $\sin(k\overline{n})=0$; if $\cos(k\overline{n})=0$, on the opposite side, the core is also apparently absent, but the sign of the field is changed.   
    \item Let us consider a \textbf{closed} waveguide $[-h,h]\times\R$, with index of refraction $n_h=\overline{n}/h$ and Dirichlet boundary conditions. By separation of variables, one can check that the Helmholtz equation can be uniquely solvable (with a choice of radiation conditions) by a function of the form $\sum_{p\in \mathbb{N}}\varphi_p(z)(a_p\cos(\sqrt{\lambda_p}x)+b_p\sin(\sqrt{\lambda_p}x))$, with $\sqrt{\lambda_p}=p\pi/h$, except if the frequency $k^*$ cancels out one of the effective frequencies $\sqrt{(n_h{k^*})^2-\lambda_p}$. In this context, these values $k^*_p=\frac{p\pi}{n_h h}=\frac{p\pi}{\overline{n}}$ for $p\in\mathbb{N}$ are called resonant frequencies.
    But these values of $k^*_p$ are exactly the roots of
    \begin{align}\label{eq:resonances}
        \sin(k\overline{n})\cos(k\overline{n})=0.
    \end{align}
    Hence, in our setting of an open waveguide with a core $[-h,h]\times\R$, with index of refraction $n_h=\overline{n}/h$, we also refer to the frequencies $k^*_p$ satisfying \eqref{eq:resonances} as \textbf{resonant frequencies}.
    \item At a resonant frequency $k^*_p$,
    we observe that the mode $v_s$ or $v_a$ corresponding to $\sqrt{\lambda}=k^*_p\overline{n}/h$ is of order 1 within the core, and of order $1/h$ outside. In some sense, such a mode only needs very little energy inside the core to transmit energy from one side of the core to the other. This is a somewhat dual effect as that of plasmonic resonances, as reported in \cite{Ebbesen}, where the amplification of the field inside a layer of resonant cavities is what allows transmission of energy from one side to the other.
    \item In Theorems \ref{thm:asym_sym} and \ref{thm:asym_anti}, in the resonant cases, we note that the first order corrections to $G$ contain an oscillating part in $z$, that resembles a guided mode. Surprisingly, this term arises as $h\to 0$ from the continuous part of the Green function.
    \item Figure \ref{fig:plotG} illustrates the dramatic change of behavior of $G$ described in Theorem \ref{thm:total_G}, as a function of the frequency $k$. In this picture, the position of the observer and the source are located on the same size of the core. The graph of $G$ essentially follows that of $2H_a$, and abruptly jumps to that of $H$ when $k$ gets close to a resonance.
    \item Refining the proofs of Theorems \ref{thm:asym_sym} and \ref{thm:asym_anti}, it can 
    be verified that the asymptotic resonant regime is still valid for frequencies  
    $k^*_p+o(h^2)$, and that the transition zone between the two regimes has a width $O(h)$ around the resonances. This is confirmed by our numerical experiments, see Figure \ref{fig:orderh}, and we take advantage of this fact for the reconstruction of $h$ in Section \ref{sec:numerical-results}, see Figure~\ref{fig:deltas}.
    \end{enumerate}
\begin{center}
\begin{figure}[ht]
\includegraphics[trim=0 0 0 0, clip, width=13cm]{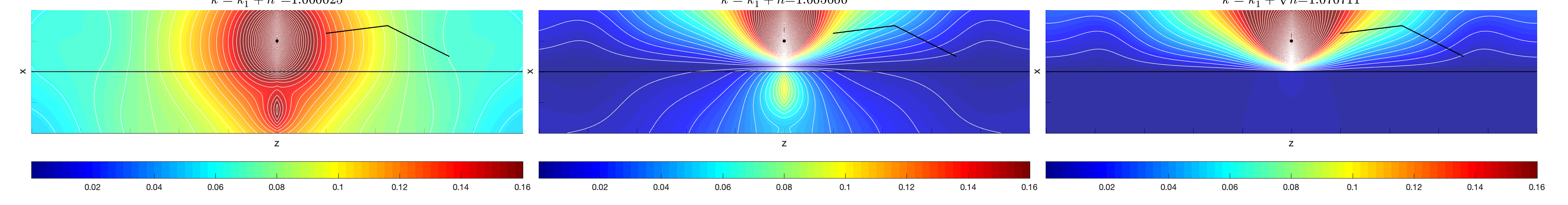}
\caption{Plots of the modulus of the field for frequencies close to the first resonance (parameters as in Figure \ref{fig:plotG}): at $k=k_1^*+h^2$ (left), at $k=k_1^*+h$ (middle), and at $k=k_1^*+\sqrt{h}$ (right). Compare the field on the left and right to the resonant and non-resonant cases in Figure \ref{fig:plotG}.}
\label{fig:orderh}
\end{figure}
\end{center}

\end{remark}

\section{Inverse problem for the detection of a high contrast core}
\label{sec:inverse_problem}

\subsection{The inverse problem}

\begin{center}
\begin{figure}[ht]
\centering
\includegraphics[trim=50 50 0 50, clip, width=10cm]{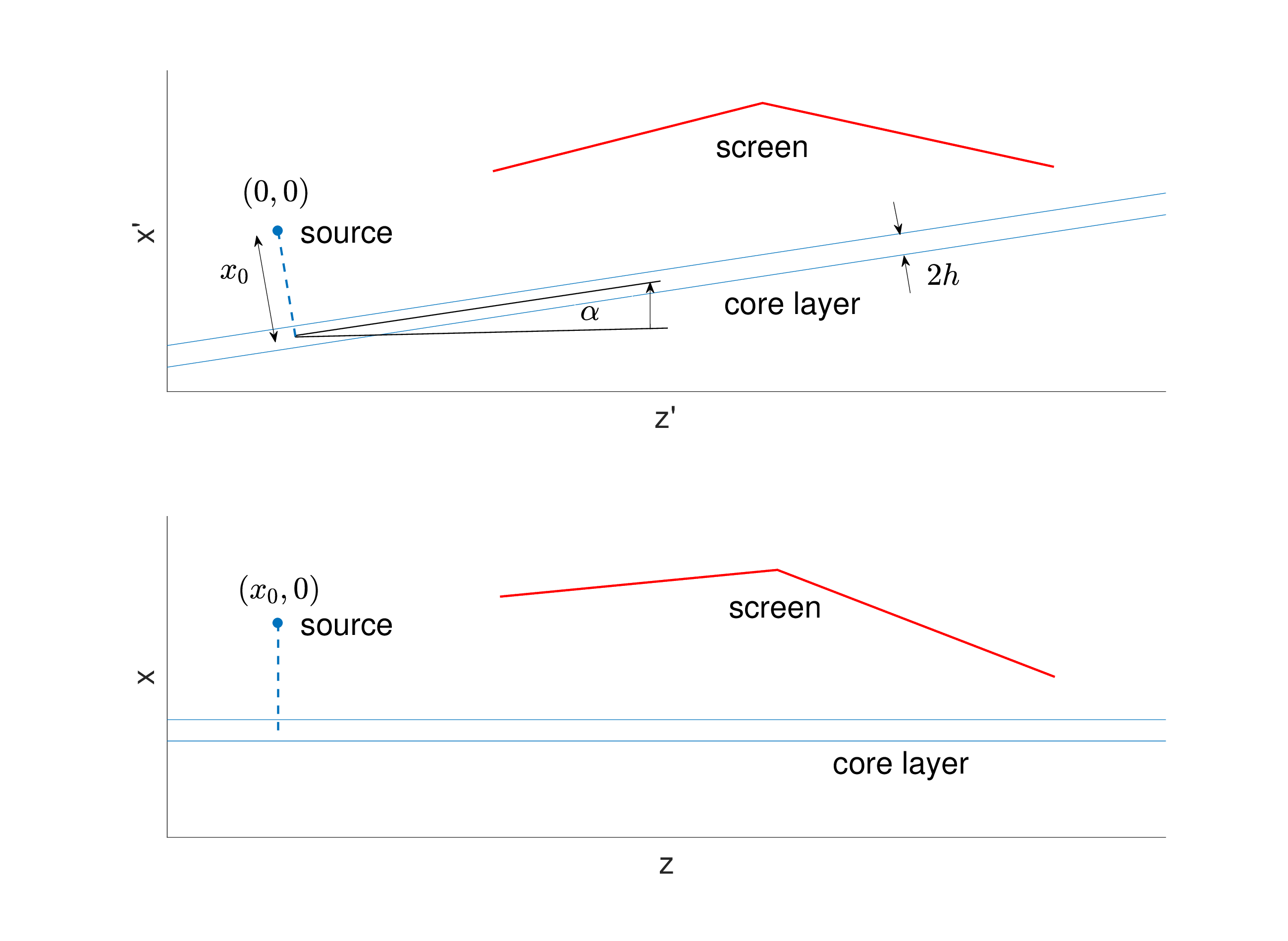}
\caption{Position of the point source, core layer and screen. 
The angle of rotation $\alpha$, the distance $x_0$ from the point source to the core (dashed line), the thickness $2h$ and the index of refraction $n_h$ of the core layer are unknowns. Measurements 
of the Green function are available on the screen $S$ for a range of frequencies $k$ (see Figure~\ref{fig:plotG}).}
\label{fig:geometry-main}
\end{figure}
\end{center}

In a reference system $Ox'z'$ (see Figure \ref{fig:geometry-main}) we consider a straight core layer of width $2h$ and optical index of refraction $n_h$, oriented along a direction that makes an angle $\alpha$ with with the $Ox'$ axis, and whose center line is at a distance $x_0$ from the origin, as shown in Figure \ref{fig:geometry-main}.
We assume that the background medium has an index $n_{cl}$ and that a source, located at the origin, emits a wave at frequency $k$, and the resulting field is measured on a screen $S\subset \R^2$. We assume that both, the source and the screen, are located strictly above the core layer. Our objective is to determine the parameters that define the core layer:
\begin{align*}
    (x_0,\alpha,h,n_h).
\end{align*}
The measurements consist in the values
\begin{align*}
    \overline{G}(x',z';0,0;k)= G(T_{x_0,\alpha}(x',z');T_{x_0,\alpha}(0,0);k)=G(T_{x_0,\alpha}(x',z');x_0,0;k),
\end{align*}
for a range of frequencies $k\in[k_{\text{min}}, k_{\text{max}}]$, that contains at least the first resonant frequency. Here, $T_{x_0,\alpha}$ denotes the affine transformation consisting in the clockwise rotation by $\alpha$, followed by the vertical translation of amplitude $x_0$, and where $G$ is the solution to \eqref{eq:Helmholtz1}-\eqref{eq:Helmholtz2} with right-hand side $f=\delta(x-x_0)\delta(z-z_0)$, and whose expression is given by \eqref{eq:def_G}. Note that, when necessary, in this section we may modify the notation to make the dependence on $k$ explicit.

Observe that the asymptotic formulas of Section \ref{sec:asymptotic} are written in the reference system $Oxz$, where the core layer is horizontal and centered about the $z$-axis, and where the point source is located at $(x_0,0)$. Hence we define
\begin{equation}\label{eq:defg}
\widetilde{S}:=T_{x_0,\alpha}(S).
\end{equation}

\subsubsection*{Steps of the identification algorithm}

\subsubsection*{$\bullet$ {Step 1: searching for the first resonance.}}
Acquire measurements $k\mapsto \overline{G}(\cdot,\cdot;0,0;k)$ on the screen,
probing the media with a sufficiently fine sample of frequencies.
According to the asymptotic results in Section \ref{sec:asymptotic}, for most frequencies, the values of the measurements are close to the values of $2H_a$. However, when $k$ is close to a resonance, the value of the measurements abruptly change to come close to the values of $H$ (see Figure~\ref{fig:plotG}). Since the values of $2H_a$ and $H$ vary smoothly in $k$ we obtain an estimate $\widehat{k}_1^*$ of $k_1^*$ as the first significant local minimum of $||\overline{G}(\cdot,\cdot;0,0;k)-H(\cdot,\cdot;0,0;k)||_{L^2(S)}$, and observe that $H$ requires only the value of $n_{cl}$ to be computed. We then estimate $\widehat{\overline{n}}=\frac{\pi}{2\widehat{k}^*_1}$.\\

\subsubsection*{$\bullet$ {Step 2: estimation of the position of the core layer.}}
Fix a non-resonant frequency $k$ (a non-integer multiple of $\widehat{k}^*_1$, found in the first step) and find the values of $x_0$ and $\alpha$ that 
minimize the gap between the measurements on the screen and  
the zero order approximation $2H_a$ given by Theorem \ref{thm:total_G}. More precisely, this amounts to finding 
\begin{equation}\label{eq:step2}
(\widehat{x}_0,\widehat{\alpha})=\arg\min_{x_0,\alpha} \|\overline{G}(x',z';0,0;k)-2H_a(T_{x_0,\alpha}(x',z');x_0,0;k)\|^2_{L^2(S)}.
\end{equation}

\subsubsection*{$\bullet$ {Step 3: estimation of the thickness and index of refraction.}}
\begin{enumerate}[(i)]
	\item One strategy is to use the first order asymptotic results of Theorem \ref{thm:total_G}, at a non-resonant frequency $k$, and estimate $h$ by
\begin{equation}\label{eq:step3}
\widehat h =\frac{1}{|\widetilde{S}|}\left|\int_{\widetilde{S}} \frac{ \overline{G}(T^{-1}_{\widehat x_0,\widehat \alpha}(x(t),z(t);0,0;k)-2H_a (x(t),z(t);\widehat x_0,0)}{(\Phi_s+\Phi_a)(x(t),z(t);\widehat x_0,0;k)}dt\right|.
\end{equation}
using the parameters $\widehat x_0$, $\widehat\alpha$ and $\widehat{\overline n}$ from the previous steps. However, our numerical experiments indicate that this approach is too sensitive to the quality of the measurements to provide a reliable approximation of $h$.
\item A second strategy consists in estimating $h$ from the width of the peaks of the graph of $\|\overline{G}\|$, around the resonant frequencies. The numerical experiments in Section \ref{sec:numerical-results} confirm that these widths are indeed proportional to $h$.
\end{enumerate}

Finally, we estimate the index of refraction in the core layer as $\widehat{n}_h=\widehat{\overline n}/\widehat{h}$.\\
We present in Section~\ref{sec:numerical-results} more details about the practical numerical implementation of the previous steps.


\subsection{Analysis of the inverse problem algorithm}

In this subsection we derive some results concerning the identification of the core layer parameters, based on the asymptotic structure of $G$ and the properties of $H$, the fundamental solution to the Helmholtz equation in the plane.

We first recall that $H(x,z;x_0,z_0)$ only depends on the distance between $(x,z)$ and $(x_0,z_0)$, so that for any isometry $T:\R^2\to\R^2$, 
\begin{align}\label{eq:isomH}
    H(x,z;x_0,z_0)=H(T(x,z);T(x_0,z_0)).
\end{align}

{\em Step 1} of the algorithm is based on the notion that the resonant and non-resonant regimes can be distinguished in the measurements and without knowing the value of any other parameter of the core layer (see Figure \ref{fig:plotG}). The lemma below shows that, asymptotically, this is indeed the case.

\begin{lemma}
    Assume that $x_0, S$ are contained in a ball of radius $R$. Then, there exists $C=C(R,k_{max})>0$ such that $\forall k\in[k_{min},k_{max}]$,
    \begin{align*}
        \inf_{(x,z)\in\widetilde{S}} |(2H_a-H)(x,z;x_0,0)|\geq C.
    \end{align*}
    If additionally, $x_0$ and $S$ are bounded away from the core layer, then we have in particular, that for $k$ non-resonant,
\begin{align*}
\lim_{h\to 0}
\inf_{(x',z')\in S}|
\overline{G}(x',z';0,0;k) - H(x',z', 0,0)|\geq C.
\end{align*}
While for $k=k^*_p$ resonant,
\eqref{eq:isomH} and of Theorem \ref{thm:total_G} imply that the measurements take the form 
\begin{align*}
    \overline{G}(x',z'; 0,0;k^*_p)&=G(T_{x_0,\alpha}(x',z');x_0,0;k^*_p)=H(T_{x_0,\alpha}(x',z');x_0,0)+O(h)\\
    &=H(x',z'; 0,0)+O(h).
\end{align*}
and therefore
\begin{align*}
\lim_{h\to 0}
\sup_{(x',z')\in S}|
\overline{G}(x',z';0,0;k^*_p) - H(x',z', 0,0)|=0.
\end{align*}
\end{lemma}

\begin{proof}The first statement follows from the fact that
\begin{align*}
(2H_a-H)(x,z,x_0,0)=-H(-x,z;x_0,0)=\frac{i}{4}\mathcal{H}\left(k n_{cl}\sqrt{(x+x_0)^2+z^2} \right)
\end{align*}
and the Hankel function never is bounded away from zero in any bounded interval (since $|\mathcal{H}(t)|$ blows up at $t=0$ and $t|\mathcal{H}(t)|^2$ is increasing for $t\in(0,\infty)$, see e.g. \cite{Watson1944}). The second part follows from Theorems \ref{thm:asym_sym} and \ref{thm:asym_anti}.
    
\end{proof}

\begin{center}
\begin{figure}[ht]
  \begin{minipage}[c]{0.4\textwidth}
\includegraphics[trim=0 15 0 25, clip, width=4.5cm]{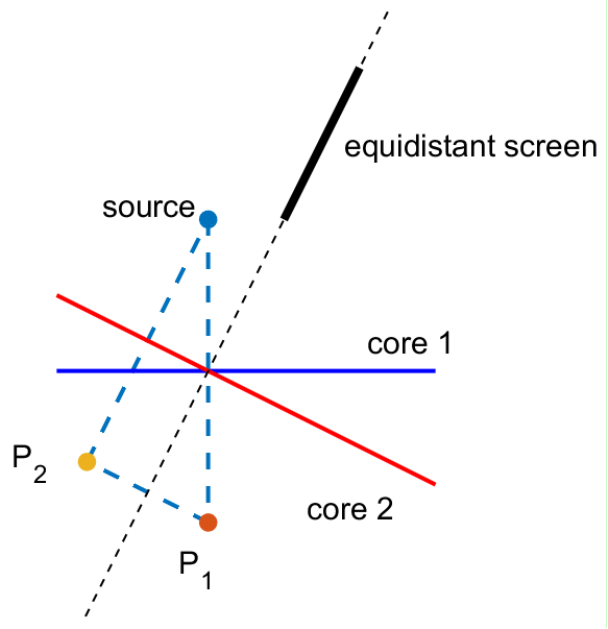}
  \end{minipage}\hfill
  \begin{minipage}[c]{0.59\textwidth}
\caption{When the reflections $P_1,P_2$ of the source with respect to the two core layers are symmetric with respect to the line containing $S$, then $H(x',z';P_1)=H(x',z';P_2)$, for all $(x',z')\in S$.}\label{fig:non-uniqueness}
  \end{minipage}
\end{figure}
\end{center}

To partially justify {\em Step 2} of the algorithm, we consider the zero order linearization in $h$ of the measurements at the non-resonant frequencies. Namely, consider two core layers with the same contrast $\overline{n}/h$ whose location is given  by the parameters $(x_1,\alpha_1)$ and $(x_2,\alpha_2)$, with $x_1,x_2>0$. Assume that $k$ is non-resonant. Then, Theorem \ref{thm:total_G} states that the respective measurements satisfy, for all $(x',z')$~in~$S$,
\begin{align*}
\overline{G}_{j}(x',z';0,0;k)=&2H_a(T_{x_j,\alpha_j}(x',z');x_j,0) + O(h)\\
=&H(T_{x_j,\alpha_j}(x',z');x_j,0)-H(T_{x_j,\alpha_j}(x',z');-x_j,0)+O(h)\\
=&H(x',z';T_{x_j,\alpha_j}^{-1}(x_j,0))-H(x',z';T_{x_j,\alpha_j}^{-1}(-x_j,0))+O(h)\\
=&H(x',z';0,0)-H(x',z';P_j)+O(h),
\end{align*}
where $P_j=T_{x_j,\alpha_j}^{-1}(-x_j,0)$ is the symmetric point to the origin with respect to the $j$-th core layer (see Figure \ref{fig:non-uniqueness}). Observe that there is a one-to-one correspondence between the coordinates of $P_j$ and the parameters $(x_j,\alpha_j)$ when $x_j>0$, thus determining $(x_j,\alpha_j)$ is equivalent to determining $P_j$.

The unique determination of the core layer parameters reduces to showing that $\overline{G}_1=\overline{G}_2$ on $S$ implies $P_1=P_2$.
Neglecting terms smaller or equal than $O(h)$, this translates into showing that 
\begin{align*}
H(x',z';P_1)=H(x',z';P_2) \text{ on } S \Rightarrow P_1=P_2.
\end{align*}
Measurements taken on a straight screen $S$ might not allow to conclude that $P_1=P_2$, as each point in $S$ could be equidistant to the points $P_1$ and $P_2$. This is exemplified in Figure \ref{fig:non-uniqueness}. On the other hand, Lemma \ref{lemma:uniquewg} below, implies that measurements on a screen with two non-parallel straight segments are enough to conclude that $P_1=P_2$.    
\begin{lemma}\label{lemma:unique1line}
Let $x_1,x_2<0$ and assume that for all $z\in[0,1],  H(0,z;x_1,z_1)-H(0,z;,x_2,z_2)=0$, then $(x_1,z_1)=(x_2,z_2)$. 
\end{lemma}
\begin{proof}
    Let $W(x,z):= H(x,z;x_1,z_1)-H(x,z;,x_2,z_2)$. Recall that  \begin{align*}
        H(0,z;x_j,z_j)=-\frac{i}{4}\mathcal{H}\left(k n_{cl}\sqrt{x_j^2+(z-z_j)^2}\right), 
    \end{align*}
    where $\mathcal{H}$ denotes the Hankel function, which is holomorphic in the complex plane cut along the negative real axis. It follows that $z\mapsto W(0,z)$ can be extended analytically to $z\in\R\times i[-\epsilon,\epsilon]$, for any $\epsilon<\min\{|x_1|,|x_2|\}$. Since $W(0,z)$ vanishes for $z\in[0,1]$, it vanishes in $\R\times i[-\epsilon,\epsilon]$, in particular
    \begin{align}\label{eq:W=0inR}
        W(0,z)=0, \forall z\in\R.
    \end{align}
    
    As $x_1,x_2<0$, we thus observe that in the upper-haf plane $P_+$, $W(x,z)$ solves the homogeneous Helmholtz equation, with standard radiation conditions, and the boundary condition \eqref{eq:W=0inR}. By uniqueness $W(x,z)=0$ in all of $P_+$.
    
    Unique continuation implies that $W(x,z)=0$ in $\R^2\setminus\{(x_1,z_1),(x_2,z_2)\}$. As each $H(x,z;x_j,z_j)$ is smooth everywhere except at $(x_j,z_j)$, this can only happen if $(x_1,z_1)=(x_2,z_2)$.
\end{proof}

\begin{lemma} \label{lemma:uniquewg}
    Assume $S\subset\R^2$ contains two non-parallel straight segments, and assume that for all $(x,z)\in S, H(x,z;x_1,z_1)-H(x,z;,x_2,z_2)=0$. Then $(x_1,x_2)=(x_2,z_2)$.
\end{lemma}

\begin{proof}
    Let $W(x,z):= H(x,z;x_1,z_1)-H(x,z;,x_2,z_2)$.
    Let $L_1,L_2$ be non-parallel lines in $\R^2$ such that $S\cap(L_1\cup L_2)$ contains two non-parallel straight segments.
    \begin{enumerate}[(i)]
        \item Case 1. Assume that $(x_1,z_1), (x_2,z_2)$ are in $ L_1$ and that $(x_1,z_1)$ is closer to $L_1\cap S$ than $(x_2,z_2)$. By analyticity of the Hankel function $\mathcal{H}$, $W(x,z)$ vanishes at all point $(x,z)$ in between $(x_1,z_1)$ and $L_1\cap S$. Since $H(x,z;x_1,z_1)$ blows up near $(x_1,z_1)$, $W(x,z)=0$ near $(x_1,z_1)$ is only possible if $(x_1,z_1)=(x_2,z_2)$.
        \item Case 2. The above argument also applies if only one of the two points belongs to either $L_1$ or $L_2$.
        \item Case 3. If neither point is in $L_1$ or $L_2$, then $W(x,z)=0$ for all $(x,z)\in S$ implies, by analyticity of $\mathcal{H}$, that $W(x,z)$ vanishes on $L_1\cup L_2$. These lines divide $\R^2$ in four regions, one of which does not contain $(x_1,z_1)$ nor $(x_2,z_2)$. We can repeat the argument in the proof of the previous Lemma to conclude that $W(x,z)=0$ in such region, and invoking again unique continuation, $W(x,z)=0$ in $\R^2\setminus\{(x_1,z_1), (x_2,z_2)\}$. The smoothness of $H(x,z; x_j,z_j)$ in all of $\R^2$, except at $(x_j,z_j)$, then implies that $(x_1,z_1)=(x_2,z_2)$.
    \end{enumerate}
\end{proof}

\section{Numerical identification of the core layer parameters}
\label{sec:numerical-results}

In this section we illustrate the numerical implementation of the identification algorithm,
proposed in Section \ref{sec:inverse_problem}, for recovering $\overline{n}, x_0,\alpha$ and $h$.

\subsection{General setting}
We consider the setting defined in Section \ref{sec:inverse_problem}. The affine linear transformation $T_{x_0,\alpha}$ introduced there is explicitly given by:
\[\left(\begin{array}{c}x\cr z\end{array}\right)=\left(\begin{array}{c}x_0\cr z_0\end{array}\right)+
\left(\begin{array}{cc}\cos\alpha&-\sin\alpha\cr \sin\alpha&\cos\alpha\end{array}\right)\left(\begin{array}{c}x'\cr z'\end{array}\right)=:T_{x_0,\alpha}(x',z').\]

For the numerical experiments, the parameters of the target core layer  are $x_0=1$ and $\alpha=\pi/20$. We choose $n_{cl}=1$, $n_h=\overline{n}/h$ with $\overline{n}=\pi/2$, so that $k_p^*=p$ for all $p\in\N$. We will run experiments for several values of $h$. In the $Oxz$ system, the screen is the union of two segments, and it is defined by 
$\tilde{S}=\{(x,z)\,\vert\, x=a_1 (z-z_m)+b_1,\ z_1<z<z_m,\ x=a_2 (z-z_m)+b_2,\ z_m<z<z_2\}$, where $z_m=(z_1+z_2)/2$, $a_1=0.1$, $b_1=0.1$, $a_2=-0.4$, $b_2=0.1$, $z_1=2$, $z_2=7$. 
Hence, the physical screen in the $Ox'z'$ reference system is $S=T_{x_0,\alpha}^{-1}(\tilde{S})$. 

To generate the synthetic data, we numerically evaluate the expression \eqref{eq:def_G} for $G(x,z;x_0,0;k)$ at the points $(x,z)\in\tilde{S}$. Hence, perfect measurements in the physical space, for each $(x',z')$ in $S$, correspond to 
\[\overline{G}(x',z';0,0;k)=G((T_{x_0,\alpha}(x',z');x_0,0;k).\] 
The continuous part of this Green function
is computed using an adapted numerical integration that uses a finer mesh of step $\Delta\tau=k\,n_{cl}/1000$ (depending on frequency $k$) near the two possible singularities $\tau=0$ and $\tau=k\,n_{cl}$ of the involved integrals. Even if in most of the considered cases, the guided part of the Green function is small compared to the continuous part (see Section \ref{sec:asymptotic}), we still compute it, in order to numerically verify the theoretical description provided in Subsection \ref{subsec:asymptotic-formulas_guided}. Finally, we add to the observations a normally distributed error, whose standard deviation is proportional to a prescribed percentage of the amplitude of the real and imaginary part of the solution $G$, at each point of the receptor. For the computation of $H$ and $2H_a$ we directly use their expressions in terms of the Hankel function, and we evaluate the Hankel function using the existing matlab implementation of it.

The identification problem consists in finding the location (angle of the waveguide core layer $\alpha$ and distance from the source to the core $x_0$), thickness $h$, 
and index of refraction $n_h$ 
of the core layer, from the previous simulated multi-frequency observations. In Figure~\ref{fig:deltas} we perform simulations for core thickness $h$ in the set $\{0.005,\,0.01,\,0.02,\,0.03,\,0.04,\,0.05\}$ and show detailed results for the cases $h=0.005$ and $h=0.05$ in Table~\ref{table:estimated-core-parameters} and Table~\ref{table:estimated-core-parameters2}.

\subsection{Step 1: searching for the first resonance}

We search for the first resonance frequency $k^*_1=1$  
by probing the medium with observations that are noisy
and measured on the receptor screen $S$, in a successively smaller and finer range of frequencies. 

\begin{center}
\begin{figure}[t]
\centering
\includegraphics[trim=80 0 120 0, clip, width=13cm, height=5.5cm]{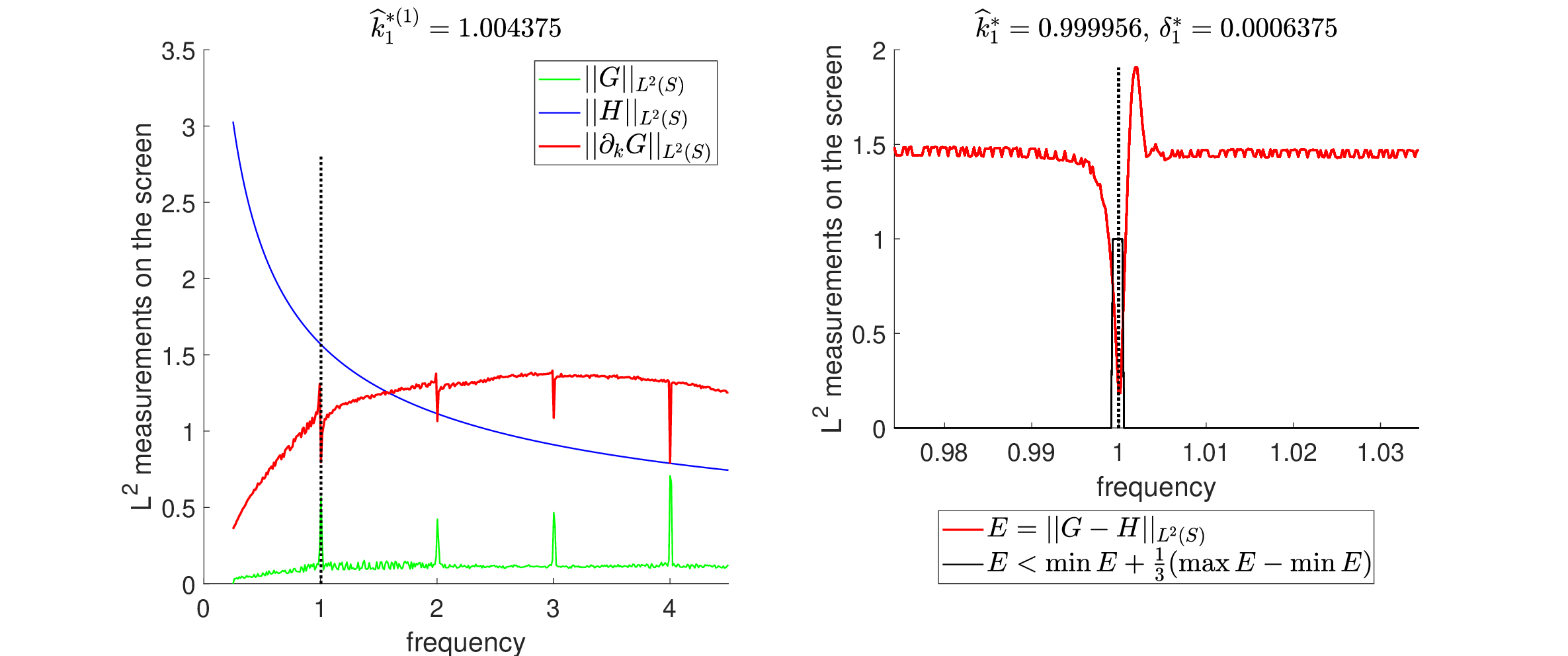}
\caption{Successive estimations $\widehat k^{*(1)}_1$ (left) and $\widehat k^*_1$ (right) of the  first radiating resonance $k^*_1$  from measurements with $5\%$ noise.  It also shows the estimated thickness $\delta^*_1$ of the first peak, which is used in Step 3 strategy (ii) of the algorithm.}
\label{fig:search-for-first-radiating-resonance}
\end{figure}
\end{center}

We start with a frequency range $k\in [k_{min},k_{max}]$, choosing $k_{min}=0.25$ and $k_{max}=4.5$ (we implicitly assume that the selected range contains the first unknown resonance frequency $k^*_1=1$ and possibly others, in this example our probing range contains the first four resonances: $k^*_p=p$, $p=1,2,3,4$). The initial frequency probing is chosen such that it does not exactly sample the resonances. And we implement the search as shown in Figure~\ref{fig:search-for-first-radiating-resonance}.

To detect the abrupt changes which occur when switching from non-resonant to resonant regime in the measurements we compute the $L^2(S)$-norm of the discrete derivative $\overline{G}(\cdot,\cdot;0,0;k+\Delta k)-\overline{G}(\cdot,\cdot;0,0;k)$.
We start the search process with a frequency step of $\Delta k=(k_{max}-k_{min})/400=0.0106$ in the whole frequency range $[k_{min},k_{max}]$. A first estimation $\widehat k^{*(1)}_1$ 
of $k^*_1$ is obtained by detecting the first local maximum of the variation prescribed above, specifically:
\[\widehat k_1^{*(1)}=\textrm{(first local)} \argmax_{k\in [k_{min},k_{max}]}  ||\overline{G}(\cdot,\cdot;0,0;k+\Delta k)-\overline{G}(\cdot,\cdot;0,0;k)||_{L^2(S)}.\]
We then refine the search in the neighborhood $K:=[0.97 k_1^{*(1)},1.03 k_1^{*(1)}]$ (in our example we refine to $\Delta k/50$), and look for the frequency that minimizes the distance between the measurements and the zero order approximation in the resonant regime (see Theorem \ref{thm:total_G}). Specifically we set:
\[\widehat k_1^*=\textrm{(first local)} \argmin_{k \in K}||\overline{G}(\cdot,\cdot;0,0;k)-H(\cdot,\cdot;0,0;k)||_{L^2(S)}.\] 
Since $k^*_1$ is the first root of $\sin(2k^*_1\overline n)=0$, we then estimate $\overline{n}$ as
$\widehat{\overline{n}}=\pi/(2\widehat k_1^*)$.

\subsection{Step 2: estimating the location of the core layer}
\begin{center}
\begin{figure}[ht]
\begin{minipage}{.63\textwidth}
\includegraphics[trim=0 0 0 0, clip, width=8cm]{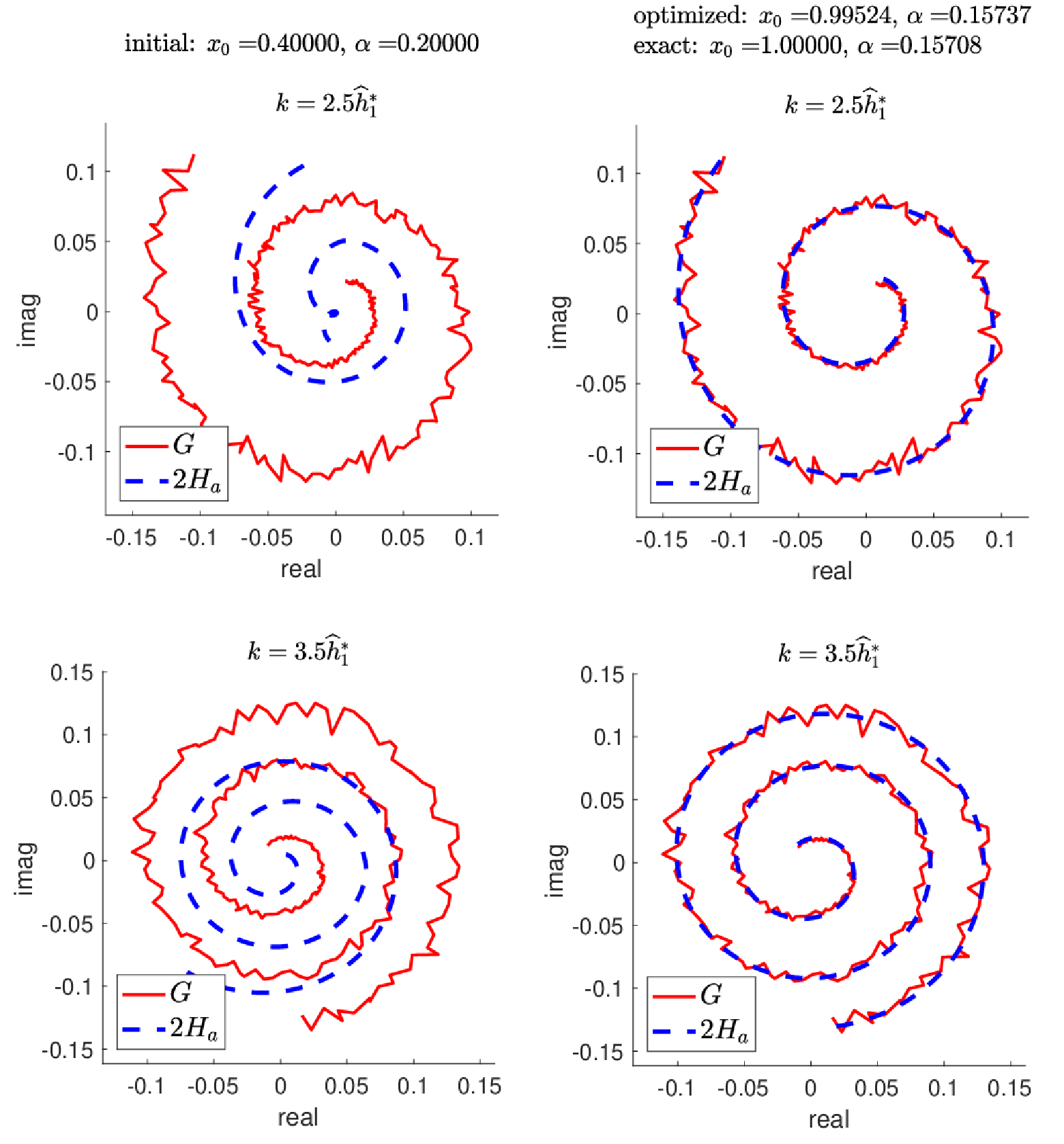}    
\end{minipage}
\begin{minipage}{.36\textwidth}
\caption{The estimation of $(x_0,\alpha)$, the core layer location parameters, is done by minimizing the difference between the measurements $G$ (real and imaginary parts) and the zero order approximation $2H_a\circ T_{x_0,\alpha}$ at the non-resonant frequencies $K=\{2.5 \widehat k_1^*,3.5\widehat k_1^*,4.5\widehat k_1^*\}$.\\
Left column: Initial guesses for $k=2.5\widehat k_1^*$ and $k=3.5\widehat k_1^* $. Right column: final approximations after optimization. Measurements with $5\%$ noise and $h=0.005$.}
\label{fig:spiral}    
\end{minipage}
\end{figure}
\end{center}

We fix a finite set $K$ of non-resonant frequencies, for our example we choose $K=\{2.5\widehat k^*_1, 3.5\widehat k^*_1, 4.5\widehat k^*_1\}$, where $\widehat k^*_1$ is the first resonant frequency estimated in step 1. 
Then we find the values of $\widehat\alpha$ and $\widehat x_0$ that minimize the $L^2(S)$ distance between the measurements and the zero order approximation in the non-resonant regime (see Theorem \ref{thm:total_G}). Specifically,
\[(\widehat x_0,\widehat\alpha)=\arg\min_{x_0,\alpha} \sum_{k_0\in K}\vert\vert \overline G(\cdot,\cdot;0,0;k_0)-2 H_a(T_{x_0,\alpha}(\cdot,\cdot);x_0,0;k_0)\vert\vert_{L^2(S)}^2.\]
We use a standard Nelder-Mead simplex algorithm \cite{Lagarias1998} called {\tt fminsearch} in the {\sc MATLAB} \cite{matlab2017a} implementation. Figure \ref{fig:spiral} illustrates how the estimation of $(\widehat x_0,\widehat\alpha)$ minimizes the distance between the measurements and $2H_a\circ T_{x_0,\alpha}$.

\subsection{Step 3: estimation of the thickness and index of refraction of the core layer}
In the last step, in order to estimate the thickness $2h$ we developed two different numerical strategies.

\begin{center}
\begin{figure}[ht]
\centering
\includegraphics[trim=30 10 30 0, clip, width=5cm]{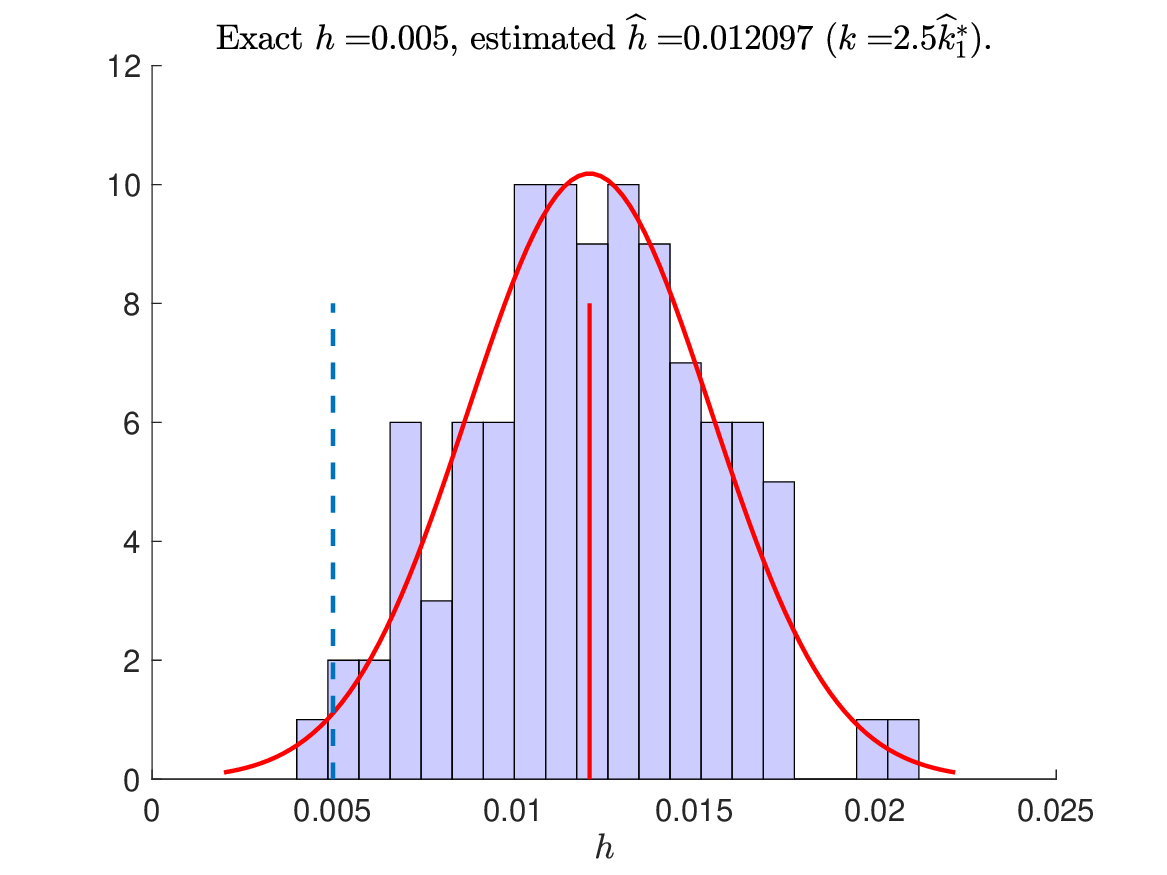}
\includegraphics[trim=30 10 30 0, clip, width=5cm]{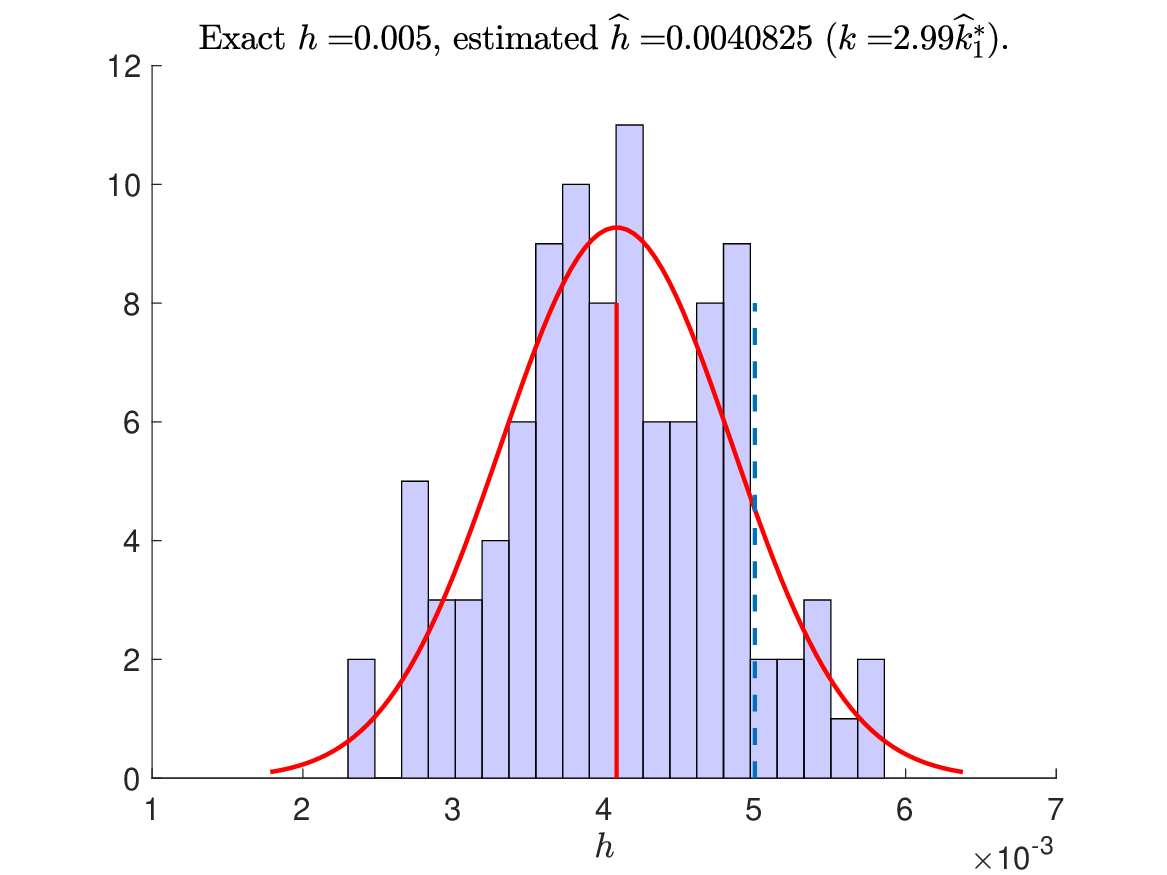}
\caption{Step 3 Strategy (i). Estimation of the thickness of the core layer by using equation \eqref{eq:step3} from multiple noise realizations, at a fixed non-resonant frequency. Examples with $k=2.5\,\widehat k^*_1$ (left) and $k=2.99\,\widehat k^*_1$ (left) ($5\%$ noise).}
\label{fig:h_estimation}
\end{figure}
\end{center}

\begin{center}
\begin{figure}[ht]
\includegraphics[trim=0 10 0 10, clip, width=13cm]{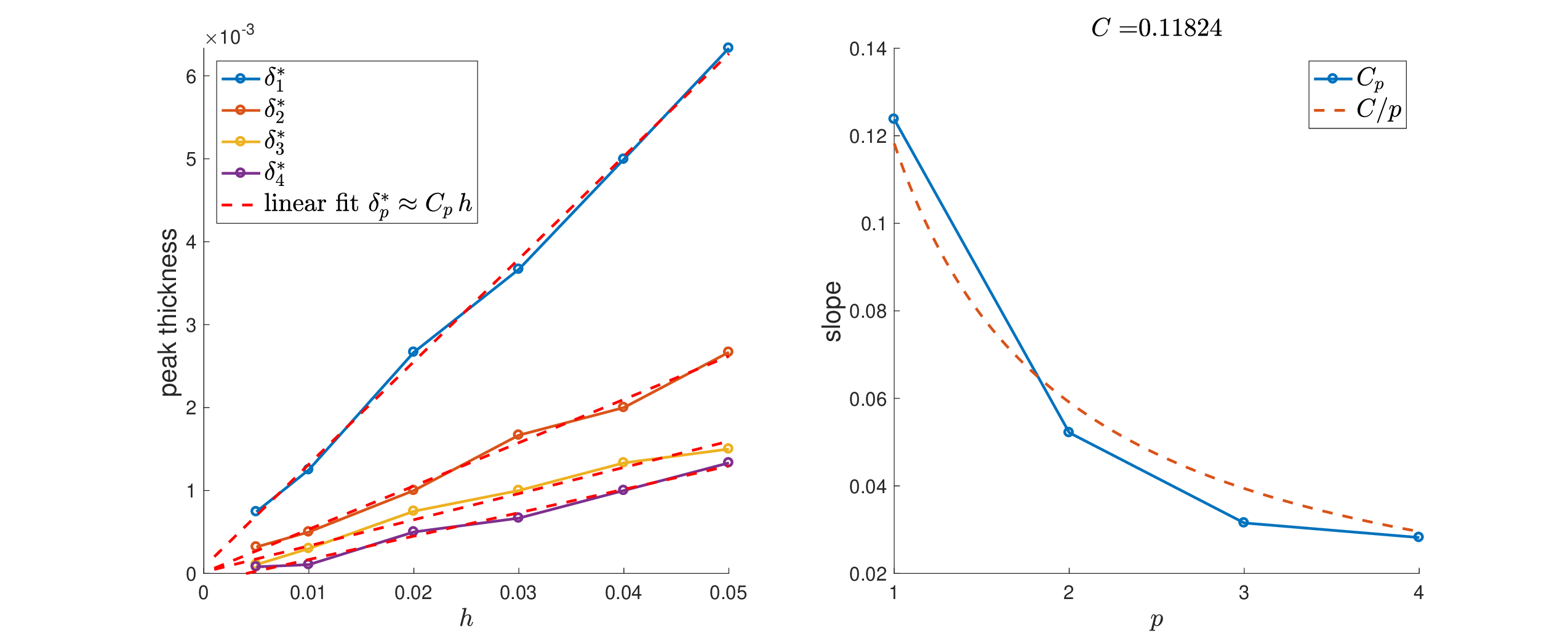}
\caption{Estimation of the thickness of the core layer in Step 3, strategy (ii). Left: Plot of $\delta^*_p(h)$ for $p=1,2,3,4$ and estimated linear fits $C_p h$. Right: 
Plot of $C_p$ for $p=1,2,3,4$ and best fit $C/p$. Data acquired with $5\%$ noise.}
\label{fig:deltas}
\end{figure}
\end{center}

{\em Strategy (i).}
We estimate $h$ from formula \eqref{eq:step3}, for a fixed non-resonant value of $k$, and several realizations of the noise. In Figure \ref{fig:h_estimation} we report the results for $k=2.5\,\widehat k^*_1$ and $k=2.99\,\widehat k^*_1$. 
We observe that the results are not very accurate, and additional numerical tests indicate that they are very sensitive to the choice of the non-resonant frequency.
Therefore, we look for an alternative way to estimate $h$ directly from the analysis of the resonance peaks of the first step.\\
{\em Strategy (ii).}  When looking for resonant frequencies in Step 1, we estimate the thickness $\delta_p^*$ of the peak corresponding to the $p$-th resonance, in the graph of the function
\[E(k)=||\overline{G}(\cdot,\cdot;0,0;k)-H(\cdot,\cdot;0,0;k)||_{L^2(S)}.\]
To this end we numerically evaluate the indicator function (see Figure~\ref{fig:search-for-first-radiating-resonance}) of the set
$B:=\left\{k:E(k)<\min_k E(k)+\beta(\max_k E(k)-\min_k E(k))\right\},$
where we arbitrarily choose $\beta=1/3$, and compute 
\[\delta_p^*=\frac{1}{2}\int_{\widehat k^*_p-\frac{1}{2}\widehat k_1^*}^{\widehat k^*_p+\frac{1}{2}\widehat k_1^*} \mathbbm{1}_B(k)dk.\]
In Figure~\ref{fig:deltas} we plot $\delta^*_p$ as a function of $h$, for $p=1,2,3,4$, and we observe a linear dependency on $h$ and an inverse dependency on $p$ of the type 
\begin{equation}\label{eq:deltas}
\delta_p^*\approx \frac{C}{p}h,
\end{equation}
with a constant $C>0$ independent of $p$ and $h$, but dependent on the choice of the set $B$ (the plot of the average slope of the curve $\delta^*_p(h)$ and the fitted $C/p$ is shown on the Figure~\ref{fig:deltas} (right)). We believe formula \eqref{eq:deltas} could be justified from the asymptotic analysis of the Green function.

In view of this, we propose to infer the value of $h$ as $h=\delta^*_1/C$.
In practice, the value of $C$ could be estimated from the $\delta^*_p$'s obtained from two controlled experiments with known values of $h$. In our experiments, we obtain $C\approx 0.11824$. We also remark that, since the computation of $E(k)$ does not require the location of the core layer, this strategy is not affected by errors in the estimation of $x_0$ and $\alpha$.\\

Finally, after estimating $h$, with strategy (i) or (ii), we estimate the index of refraction by 
\[\widehat n_h=\overline n/\widehat h.\]

A summary of the numerical results of the whole algorithm is presented in Table~\ref{table:estimated-core-parameters} for $h=0.005$ and Table~\ref{table:estimated-core-parameters2} for $h=0.05$. In both Tables, Step 3 strategy (i) was implemented with $k=2.99 \widehat k^*_1$.

These results indicate that the values of $k_1^*$ and $\overline n$ can be recovered with high accuracy even with $8\%$ noise in the measurements.
The results also show an adequate reconstruction of the location parameters of the core layer, though these are more sensitive to the noise level.
Recovery of the thickness parameter is much harder. In particular, when $h$ is smaller than the noise level, as in the case of Table~\ref{table:estimated-core-parameters}, strategy (i) is doomed to fail, because the first order term in the asymptotic formula from Theorem \eqref{thm:total_G}, from which equation \eqref{eq:step3} is derived is drowned in the noise. In view of this, the results obtained with strategy (ii) are quite satisfying and show reasonable robustness to noise.

\begin{table}
{\small
\begin{tabular}{c|c|c|c|c|c|c|c|c}
\toprule
&\multicolumn{2}{c|}{1st step}&\multicolumn{2}{c}{2nd step}&\multicolumn{2}{c|}{3th step (i)}&\multicolumn{2}{c}{3th step (ii)}\cr
\hline
parameter\!\!&\!\!\! $k^*_1$ \!\!&\!\!\! $\overline n=n_h\,h$ \!\!&\!\!\! $x_0$\!\!&\!\!\! $\alpha$\!\!&\!\!\! $h$ \!\!&\!\!\! $n_h$ \!\!&\!\!\! $h$\!\!&\!\!\! $n_h$\cr
\hline
exact value \!\!&\!\!\! 1\!\!&\!\!\! 1.5708 \!\!&\!\!\!1 \!\!&\!\!\! 0.15708\!\!&\!\!\! 5$\times 10^{-3}$\!\!&\!\!\! 314.2 \!\!&\!\!\!5$\times 10^{-3}$\!\!&\!\!\! 314.2\cr
\hline
\!\!&\!\!\!\!\!&\!\!\!\!\!&\!\!\!\!\!&\!\!\!\!\!&\!\!\! $\times 10^{-3}$\!\!&\!\!\! \!\!&\!\!\!$\times 10^{-3}$ \!\!&\!\!\!\cr
3\% noise\!\!&\!\!\! 0.99996\!\!&\!\!\! 1.5709 \!\!&\!\!\! 0.99710\!\!&\!\!\! 0.16068\!\!&\!\!\! 4.0810\!\!&\!\!\! 384.94 \!\!&\!\!\! 5.3920\!\!&\!\!\! 291.35\cr
{\it (rel. error)}\!\!&\!\!\! {\it (0.004\%)}\!\!&\!\!\! {\it (0.006\%)} \!\!&\!\!\! {\it (0.2\%)} \!\!&\!\!\! {\it (2.3\%)}\!\!&\!\!\! {\it (18\%)}\!\!&\!\!\! {\it (23\%)} \!\!&\!\!\! {\it (8\%)} \!\!&\!\!\! {\it (7\%)}\cr
5\% noise\!\!&\!\!\! 0.99996\!\!&\!\!\! 1.5709 \!\!&\!\!\! 1.0164\!\!&\!\!\! 0.17012\!\!&\!\!\! 2.8510\!\!&\!\!\! 550.91 \!\!&\!\!\! 6.2900\!\!&\!\!\! 249.73\cr
{\it }\!\!&\!\!\! {\it (0.004\%)}\!\!&\!\!\! {\it (0.006\%)} \!\!&\!\!\! {\it (1.6\%)} \!\!&\!\!\! {\it (8\%)}\!\!&\!\!\! {\it (43\%)}\!\!&\!\!\! {\it (75\%)} \!\!&\!\!\! {\it (25\%)} \!\!&\!\!\! {\it (21\%)}\cr
8\% noise\!\!&\!\!\! 0.99995\!\!&\!\!\! 1.5709 \!\!&\!\!\!0.99015\!\!&\!\!\! 0.15117\!\!&\!\!\! 1.1566\!\!&\!\!\! 135.79 \!\!&\!\!\!6.2900\!\!&\!\!\! 249.68\cr
{\it }\!\!&\!\!\! {\it (0.004\%)}\!\!&\!\!\! {\it (0.006\%)} \!\!&\!\!\! {\it (1.0\%)} \!\!&\!\!\! {\it (3.8\%)}\!\!&\!\!\! {\it (131\%)}\!\!&\!\!\! {\it (56\%)} \!\!&\!\!\! {\it (25\%)} \!\!&\!\!\! {\it (20\%)}\cr
\bottomrule
\end{tabular}
}

\caption{Results of the recovering algorithm presented in Section~\ref{sec:numerical-results} for the case  $h=0.005$ with $5$ significative digits. Exact and estimated core layer parameters using the three-step algorithm proposed in this section. The first resonant frequency $k^*_1$ and the scaled index of refraction $\overline n=n_h\, h$ in the first step, the distance from the source $x_0$ and inclination angle $\alpha$ in the second step, and the thickness $h$ and 
index of refraction $n_h$ in the third step for two different strategies. }
\label{table:estimated-core-parameters}
\end{table}

\begin{table}
{\small
\begin{tabular}{c|c|c|c|c|c|c|c|c}
\toprule
&\multicolumn{2}{c|}{1st step}&\multicolumn{2}{c}{2nd step}&\multicolumn{2}{c|}{3th step (i)}&\multicolumn{2}{c}{3th step (ii)}\cr
\hline
parameter\!\!&\!\!\! $k^*_1$ \!\!&\!\!\! $\overline n=n_h\,h$ \!\!&\!\!\! $x_0$\!\!&\!\!\! $\alpha$\!\!&\!\!\! $h$ \!\!&\!\!\! $n_h$ \!\!&\!\!\! $h$\!\!&\!\!\! $n_h$\cr
\hline
exact value \!\!&\!\!\! 1\!\!&\!\!\! 1.5708 \!\!&\!\!\!1 \!\!&\!\!\! 0.15708\!\!&\!\!\! 5$\times 10^{-2}$\!\!&\!\!\! 31.420 \!\!&\!\!\!5$\times 10^{-2}$\!\!&\!\!\! 31.420\cr
\hline
\!\!&\!\!\!\!\!&\!\!\!\!\!&\!\!\!\!\!&\!\!\! \!\!&\!\!\! $\times 10^{-2}$\!\!&\!\!\! \!\!&\!\!\!$\times 10^{-2}$ \!\!&\!\!\!\cr
3\% noise\!\!&\!\!\! 1.0012\!\!&\!\!\! 1.5690 \!\!&\!\!\! 0.97271\!\!&\!\!\! 0.16495\!\!&\!\!\! 4.3365\!\!&\!\!\! 36.181 \!\!&\!\!\! 5.0921\!\!&\!\!\! 30.812\cr
{\it (rel. error)}\!\!&\!\!\! {\it (0.1\%)}\!\!&\!\!\! {\it (0.0003\%)} \!\!&\!\!\! {\it (3\%)} \!\!&\!\!\! {\it (5\%)}\!\!&\!\!\! {\it (13\%)}\!\!&\!\!\! {\it (15\%)} \!\!&\!\!\! {\it (2\%)} \!\!&\!\!\! {\it (2\%)}\cr
5\% noise\!\!&\!\!\! 1.0012\!\!&\!\!\! 1.5690 \!\!&\!\!\! 0.78061\!\!&\!\!\! 0.10522\!\!&\!\!\! 4.7046\!\!&\!\!\! 33.350 \!\!&\!\!\! 5.0921\!\!&\!\!\! 30.812\cr
{\it }\!\!&\!\!\! {\it (0.1\%)}\!\!&\!\!\! {\it (0.0003\%)} \!\!&\!\!\! {\it (22\%)} \!\!&\!\!\! {\it (33\%)}\!\!&\!\!\! {\it (6\%)}\!\!&\!\!\! {\it (6\%)} \!\!&\!\!\! {\it (2\%)} \!\!&\!\!\! {\it (2\%)}\cr
8\% noise\!\!&\!\!\! 1.0012\!\!&\!\!\! 1.5690 \!\!&\!\!\!0.75190\!\!&\!\!\! 0.064378\!\!&\!\!\! 6.2134\!\!&\!\!\! 25.251 \!\!&\!\!\!5.3916\!\!&\!\!\! 29.100\cr
{\it }\!\!&\!\!\! {\it (0.1\%)}\!\!&\!\!\! {\it (0.0003\%)} \!\!&\!\!\! {\it (25\%)} \!\!&\!\!\! {\it (59\%)}\!\!&\!\!\! {\it (24\%)}\!\!&\!\!\! {\it (20\%)} \!\!&\!\!\! {\it (8\%)} \!\!&\!\!\! {\it (7\%)}\cr
\bottomrule
\end{tabular}
}
\caption{Results of the recovering algorithm presented in Section~\ref{sec:numerical-results} for the case  $h=0.05$. Same description as in Table~\ref{table:estimated-core-parameters}. }
\label{table:estimated-core-parameters2}
\end{table}

\section{Proofs of the asymptotic formulas}
\label{sec:asymptotic-formulas}

\subsection{Proof of Theorem \ref{thm:asym_guided}: Asymptotic analysis of the guided component.}\label{subsec:asymptotic-formulas_guided}

We show that the guided modes vanish quickly as we move away from the core layer, and do not contribute to the main terms of the asymptotic expansion.

\begin{proof}[Proof of Theorem \ref{thm:asym_guided}]
To analyze what happens to the guided modes let us first study the equation for $\{\lambda_{s,j}\}_1^{J_s}$, the roots associated to the symmetric guided modes. Namely, let us study the equation 
\begin{align}
\label{eq:lambdaroots_s}&\sqrt{d^2-\lambda}-\sqrt{\lambda}\tan(h\sqrt{\lambda})=0&\wedge& \quad 0<\lambda < d^2,\\ 
\nonumber \Leftrightarrow \quad &d^2-\lambda=\lambda\tan^2(h\sqrt{\lambda}) &\wedge& \quad \tan(h\sqrt{\lambda})>0,\\
 \nonumber\Leftrightarrow \quad &h^2d^2=h^2\lambda(1+\tan^2(\sqrt{h^2\lambda}))=h^2\lambda\sec^2(\sqrt{h^2\lambda})&\wedge& \quad \tan(\sqrt{h^2\lambda})>0.
\end{align}
By letting $y^2=h^2\lambda$ and recalling that $h^2d^2=h^2k^2(n_h^2-n_{cl}^2)=k^2\overline n^2-h^2k^2n_{cl}^2$, we obtain the equivalent equation on $y$,
\begin{align}\label{eq:yroots_s}
y^2\sec^2(y)=k^2\overline n^2-h^2k^2n_{cl}^2 ~ \wedge ~ \tan(y)>0, y>0,
\end{align}
and we denote by $\{y_{s,j}(h)\}_{j=1}^{J_s}$ its roots. Define $L^2:=k^2\overline{n}-h^2k^2n_{cl}^2$ and let us compute $p=J_s$, the number of roots of \eqref{eq:yroots_s}. 

\begin{figure}[ht]
\begin{minipage}{.44\textwidth}
	\centering
	\includegraphics[trim=120 50 40 60,scale=0.25]{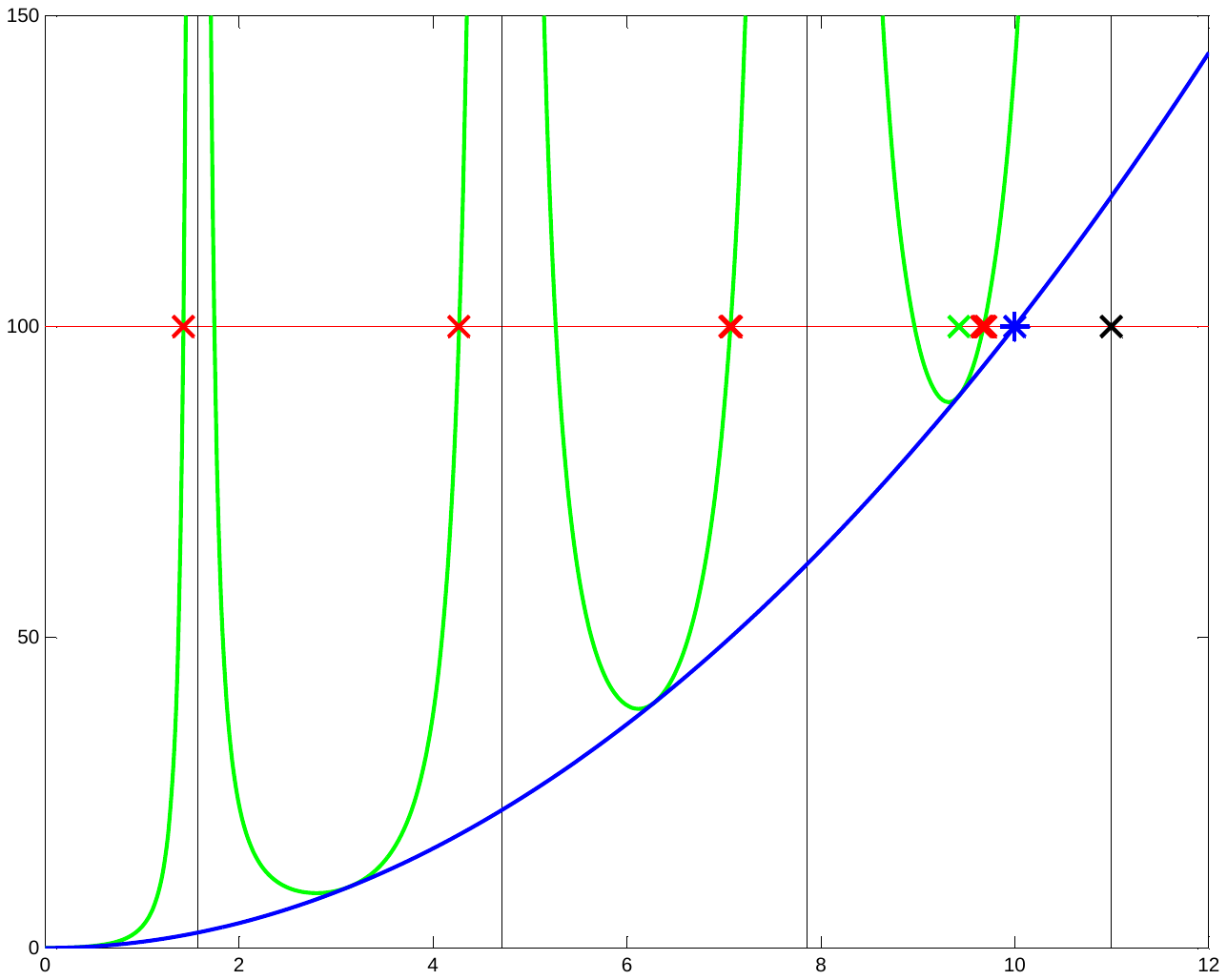}
\end{minipage}
\begin{minipage}{.55\textwidth}
    \caption{Graphs of $y^2$ (blue line) and $y^2\sec^2(y)$ (green line). Given a height $L^2$ (red line) the 
red X's mark the roots of $y^2\sec^2(y)=L^2$ at which $\tan(y)> 0$, the blue * marks the root of
$y^2=L^2$ and the black X marks $y=(2p-1)\pi/2$ for $p$ equal to the number of red roots.
The last red root is in the interval $((2p-3)\pi/2,(2p-1)\pi/2)$ if and only if $(p-1)\pi< L \leq p\pi$,
i.e. $p-1< L/\pi\leq p$ and $p=\ceil{L/\pi}$.}\label{fig:roots_s}
\end{minipage}
\end{figure}

Since $\sec^2(y)\geq 1$, and $\sec^2(y)=1$ only when $y$ is an integer multiple of $\pi$, the largest root $y_{s,p}(h)$ is related to the root of $y^2=L^2$ through the relation: $y_{s,p}(h)
\in (p\pi-3\pi/2,p\pi-\pi/2)$ if and only if $(p-1)\pi<L\leq p\pi$ (see Figure \ref{fig:roots_s}). Therefore $(p-1)<L/\pi \leq p$, and $p=\ceil{L/\pi}$, where $\ceil{x}=\min\{n\in\N:n\geq x\}$.
We conclude that equation \eqref{eq:yroots_s} admits exactly $J_s=\ceil[\Big]{\frac{1}{\pi}\sqrt{k^2\overline n^2-h^2k^2n^2_{cl}}}$ roots, and since $t\mapsto \ceil{t}$ is piece-wise constant and left continuous, for $h>0$ small, $J_s=\ceil[\Big]{\frac{1}{\pi}k\overline n}$.

If we let $y_{s,j}:=y_{s,j}(0)$, equation \eqref{eq:yroots_s} also implies, for all $j\in\{1,...,J_s\}$,
\begin{align}\label{eq:symyj<knbar}
y_{s,j}(h)=y_{s,j}+O(h^2), \quad 
|\cos(y_{s,j})|<1, \quad  0<y_{s,j}<k\overline n.
\end{align}

Recapitulating, the roots $\{\lambda_{s,j}\}_1^{J_s}$ of equation \eqref{eq:lambdaroots_s} satisfy: $J_s = \ceil{k\overline n/\pi}$, and $\lambda_{s,j}=y^2_{s,j}(h)/h^2=y^2_{s,j}/h^2+O(1)$ for $h>0$ small.

For the antisymmetric guided modes the situation is similar. The equation for the roots $\{\lambda_{a,j}\}_1^{J_a}$ is,
\begin{align*}
&\sqrt{d^2-\lambda}+\sqrt{\lambda}\cot(h\sqrt{\lambda})=0&\wedge& \quad 0<\lambda < d^2,\\
\Leftrightarrow \quad &h^2d^2-h^2\lambda=h^2\lambda\cot^2(\sqrt{h^2\lambda})&\wedge& \quad \cot(\sqrt{h^2\lambda})<0,
\end{align*}
letting $y^2=h^2\lambda$ we obtain the related equation,
\begin{align*}
\quad &k^2\overline n^2-h^2k^2n_{cl}^2=y^2(1+\cot^2(y))=\left(\frac{y}{\sin(y)}\right)^2\quad \quad \wedge \quad \cot(y)<0, y>0.
\end{align*}

\begin{figure}[ht]
\begin{minipage}{.44\textwidth}
	\centering
	\includegraphics[trim=120 50 40 70,scale=0.25]{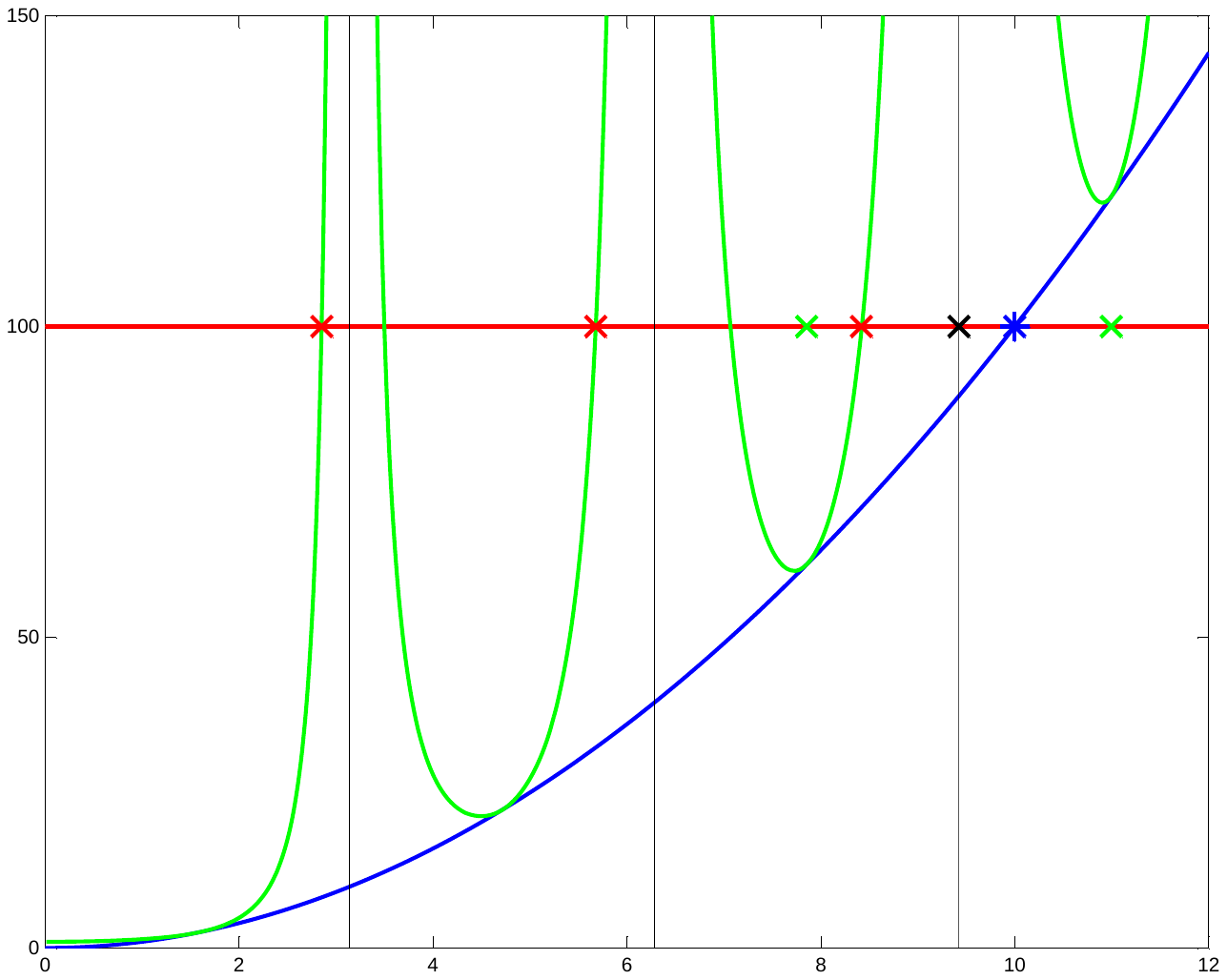}
\end{minipage}
\begin{minipage}{.55\textwidth}
    \caption{Graphs of $y^2$ (blue line) and $(y/\sin(y))^2$ (green line). Given a hight $L^2$ (red line) the 
red X's mark the roots of $(y/\sin(y))^2=L^2$ at which $\cot(y)< 0$, the blue * marks the root of
$y^2=L^2$ and the black X marks $y=p\pi$ for $p$ equal to the number of red roots.
The last red root is in the interval $((p-1)\pi,p\pi)$ if and only if $(p-1/2)\pi< L \leq (p+1/2)\pi$,
i.e. $p-1< L/\pi-1/2\leq p$ and $p=\ceil{L/\pi-1/2}$.}\label{fig:roots_a}
\end{minipage}
\end{figure}

For this equation there are exactly $J_a=\ceil[\Big]{\frac{1}{\pi}\sqrt{k^2\overline n^2-h^2k^2n^2_{cl}}-1/2}$ roots (see Figure \ref{fig:roots_a}), and for $h>0$ small  $J_a=\ceil[\Big]{\frac{1}{\pi}k\overline n-1/2}$.
If $\{y_{a,j}(h)\}_{j=1}^{J_a}$ are the roots of
\begin{align*}
&\left(\frac{y}{\sin(y)}\right)^2=k^2\overline n^2-h^2k^2n_{cl}^2\quad \wedge \quad \cot(y)<0, y>0,
\end{align*}
letting $y_{a,j}:=y_{a,j}(0)$, then $y_{a,j}(h)=y_{a,j}+O(h^2)$ for $h\sim0$. Also
\begin{align}
    |\sin(y_{a,j})|<1 \textnormal{ and } 0<y_{a,j}<k\overline{n}, \forall j=1,...,J_a.
\end{align}

For the roots $\{\lambda_{a,j}\}_1^{J_a}$ of the original equation, 
this means that there are $J_a=\ceil{k\overline n/\pi-1/2}$ roots 
and they satisfy $\lambda_{a,s}=y^2_{a,s}(h)/h^2=y^2_{a,s}/h^2+O(1)$.\\

Now let us estimate the elements in the guided component of the Green's function.
For $c>0$ and $|x_0|\geq c$, if $0<h<c/2$ is small, the behavior of 
the $\lambda_{m,j}$ imply that
$\sqrt{d^2-\lambda_{m,j}}=\frac{1}{h}(1+O(h^2))\sqrt{k^2\overline n^2-y^2_{m,j}}$, therefore
\begin{align*}
v_s(x_0,\lambda_{s,j}) &= \cos(h\sqrt{\lambda_{s,j}})e^{-\sqrt{d^2-\lambda_{s,j}}(|x_0|-h)}\\
&=(1+O(h^2))\cos(y_{s,j})e^{-\frac{|x_0|}{h}(1+O(h))\sqrt{k^2\overline n^2-y^2_{j,s}}},\\
v_a(x_0,\lambda_{a,j}) &= \sin(h\sqrt{\lambda_{a,j}})e^{-\sqrt{d^2-\lambda_{a,j}}(|x_0|-h)}\\
&=(1+O(h^2))\sin(y_{a,j})e^{-\frac{|x_0|}{h}(1+O(h))\sqrt{k^2\overline n^2-y^2_{a,s}}}.
\end{align*}
Also,  $k\beta_{m,j}=\sqrt{k^2n_h^2-\lambda_{m,j}}=\frac{1}{h}(1+O(h^2))\sqrt{k^2\overline n^2-y^2_{m,j}}$, hence
$\frac{\sqrt{d^2-\lambda_{m,j}}}{2ik\beta_{m,j}}=1+O(h^2).$
Additionally,
\begin{align*}
v_m(x,\lambda_{m,j})=O(1), \quad e^{ik\beta_{m,j}|z-z_0|}=O(1),
&\quad \text{ and }\quad  1+h\sqrt{d^2-\lambda_{m,j}}=O(1).
\end{align*}
All the estimates above are uniformly bounded for $x,z\in\R$.
Replacing such estimates in the expression of the guided component of the Green's function, we obtain 
\begin{align*}
G_{m,g}(x,z;x_0,z_0)=O(1)\sum_{j=1}^{J_m} e^{\frac{-|x_0|}{h}(1+O(h))
\sqrt{k^2\overline n^2-y^2_{m,j}}}.
\end{align*}
Since $\min_j\{k^2\overline n^2-y^2_{m,j}\}>0$ only depends on $k$ and $\overline{n}$, and the $O(1), O(h)$ estimates above are uniform on $x,z\in\R$, we obtain the desired inequality~\eqref{eq:estim_guided}.
\end{proof}



\subsection{Proof of Theorem \ref{thm:asym_sym}: Asymptotic analysis of the continuous symmetric component.}

We are only interested in an asymptotic formula of the Green function when $x$ and $x_0$ are outside the core of the waveguide. To obtain an asymptotic expression we will have to separate the cases when $\sin(2k\overline n)\neq0$ and $\sin(2k\overline n)=0$, and the expression may depend on $x,x_0$ being on the same or opposite sides of the core. We restrict our analysis to $z$ bounded away from $z_0$. We believe that it should be possible to replace this technical condition by only requiring that $(x,z)$ be bounded away from the source $(x_0,z_0)$.

Recall that we defined $\Phi_s$ and $\Psi_s$ in Theorem \ref{thm:asym_sym}. The statement (\ref{enum:sym1}) of Theorem \ref{thm:asym_sym} corresponds to the  following lemma.

\begin{lemma}\label{lemma:csym_nr}
Fix $C>c>0$ and assume that $c\leq x_0\leq C$. Also assume that $\sin(k\overline n)\neq0$, then there exists $\kappa$ only depending on $k,\overline{n}, c, C$ and $n_{cl}$ such that 
\begin{align*}
\begin{split}
&\| (G_{s,c}-[H_a+h\Phi_s])(\cdot,\cdot;x_0,z_0)\|_{L^\infty(\{c\leq x\leq C\}\times \{c\leq |z-z_0|\})}\leq \kappa h^2,\\
&\text{ and }\\
&\| (G_{s,c}-[-H_a+h\Phi_s])(\cdot,\cdot;x_0,z_0)\|_{L^\infty(\{-C\leq x\leq -c\}\times \{c\leq |z-z_0|\})}\leq \kappa h^2.
\end{split}
\end{align*}
\end{lemma}

\begin{proof}
The continuous symmetric part has the expression, 
\begin{align*}
G_{s,c}(x,z;x_0,z_0)=&\frac{1}{2\pi}\int_{d^2}^\infty v_s(x,\lambda)v_s(x_0,\lambda)
\frac{e^{ik\beta|z-z_0|}}{2ik\beta}\frac{\sqrt{\lambda-d^2}}{(\lambda-d^2)+d^2 \sin^2(h\sqrt{\lambda})} d\lambda
\end{align*}
and considering the change of variable $k^2\tau^2=\lambda-d^2=\lambda-k^2n_h^2+k^2n^2_{cl}=k^2n_{cl}^2-k^2\beta^2$, then
(multiplying and dividing by $h^2$),
\begin{align}\label{eq:Gsymm_form_tauh}
\begin{split}
G_{s,c}(x,z;x_0,z_0)=&\frac{1}{2\pi}\int_{0}^\infty hv_s(x,k^2\tau^2+d^2)hv_s(x_0,k^2\tau^2+d^2)\cdot\\
&\cdot\frac{e^{ik\sqrt{n^2_{cl}-\tau^2}|z-z_0|}}{2ik\sqrt{n^2_{cl}-\tau^2}}
\frac{k\tau}{h^2k^2\tau^2+h^2d^2 \sin^2(h\sqrt{k^2\tau^2+d^2})} 2k^2\tau d\tau
\end{split}
\end{align}
and for $h<c\leq |x|\leq C$, since $Q=\sqrt{\lambda-d^2}=k\tau$, we have that
\begin{align*}
hv_s(x,k^2\tau^2+d^2)&=
h\cos(h\sqrt{k^2\tau^2+d^2})\cos(k\tau(|x|-h))\\
&-h\sqrt{k^2\tau^2+d^2}\sin(h\sqrt{k^2\tau^2+d^2})\frac{\sin(k\tau(|x|-h))}{k\tau}.
\end{align*}
We will estimate this quantity by abusing some notation, we will denote by $O(h^m)P(\tau)$ any function that, for $h>0$ small, is bounded by $h^m$ times a polynomial in $\tau$ (uniformly bounded with respect to $|x|,|x_0|$ in $[c,C]$).
Then,
\begin{align*}
h\cos(k\tau(|x|-h))&= h\cos(k\tau|x|)\cos(k\tau h)+h\sin(k\tau|x|)\sin(k\tau h)\\
&=h\cos(k\tau|x|)+O(h^2)P(\tau),\\
\sin(k\tau(|x|-h))&= \sin(k\tau|x|)\cos(k\tau h)-\cos(k\tau|x|)\sin(k\tau h)\\
&= \sin(k\tau|x|)(1+O(h^2)P(\tau))-
k\tau h\cos(k\tau|x|)(1+O(h^2)P(\tau))\\
&= \sin(k\tau|x|)-k\tau h\cos(k\tau|x|) +O(h^2)P(\tau).
\end{align*}
Since $h^2d^2=k^2\overline n^2+O(h^2)$, then 
\begin{align*}
h\sqrt{k^2\tau^2+d^2}&= k\overline n +O(h^2)P(\tau),\\
\cos(h\sqrt{k^2\tau^2+d^2})&=\cos(k\overline n)+O(h^2)P(\tau),\\
\sin(h\sqrt{k^2\tau^2+d^2})&=\sin(k\overline n)+O(h^2)P(\tau).
\end{align*}
All of this together implies that 
\begin{align*}
hv_s(x,k^2\tau^2+d^2)&=-k\overline n\sin(k\overline n)\frac{\sin(k\tau |x|)}{k\tau}\\
&+h\big(\cos(k\overline n)+k\overline n\sin(k\overline n)\big)\cos(k\tau |x|)+O(h^2)P(\tau).
\end{align*}
And a similar estimate holds for $v_s(x_0,k^2\tau^2+d^2)$.

Additionally, since  $\sin(k\overline{n})\neq 0$, then for all $h>0$ sufficiently small,  $(h^2k^2\tau^2+h^2d^2 \sin^2(h\sqrt{k^2\tau^2+d^2})$ is bounded away from zero uniformly in $\tau$, and we can obtain that 
\begin{align*}
&\frac{2k^3\tau^2}{h^2k^2\tau^2+h^2d^2 \sin^2(h\sqrt{k^2\tau^2+d^2})}-\frac{2k^3\tau^2}{k^2\overline n^2\sin^2(k\overline n)}= O(h^2)P(\tau),
\end{align*}
Replacing this estimates in equation \eqref{eq:Gsymm_form_tauh} and recalling that the $v_s$ are uniformly bounded for $|x|,|x_0|$ in $[c,C]$, we obtain
\begin{align*}
G_{s,c}(x,z;x_0,z_0)=&\frac{1}{2\pi}\int_{0}^\infty
\sin(k\tau|x|)\sin(k\tau|x_0|)
\frac{e^{ik\sqrt{n^2_{cl}-\tau^2}|z-z_0|}}{i\sqrt{n^2_{cl}-\tau^2}}d\tau\\
&+ h \Phi_s(x,z;x_0,z_0)+
\int_{0}^\infty \frac{e^{ik\sqrt{n^2_{cl}-\tau^2}|z-z_0|}}{i\sqrt{n^2_{cl}-\tau^2}}O(h^2)P(\tau)d\tau,
\end{align*}
where $O(h^2)P(\tau)$ can be bounded, for $h>0$ small, by $h^2$ times a polynomial on $\tau$, uniformly for $|x|,|x_0|$ in $[c,C]$ and in $z, z_0$ in $\R$. We also assumed $|z-z_0|\geq c$, which allow us to compute the integral 
\begin{align*}
\int_{0}^\infty \frac{e^{ik\sqrt{n^2_{cl}-\tau^2}|z-z_0|}}{i\sqrt{n^2_{cl}-\tau^2}}O(h^2)P(\tau)d\tau=O(h^2),
\end{align*}
uniformly for $|x|,x_0\in [c,C]$ and for $|z-z_0|\geq c$. Finally, using Lemma  \ref{lemma:H} we conclude the proof.
\end{proof}

The statement (\ref{enum:sym2}) of Theorem \ref{thm:asym_sym} corresponds to the  following lemma.

\begin{lemma}\label{lemma:csym_r}
Fix $C>c>0$ and assume that $c\leq x_0\leq C$. Also assume that $\sin(k\overline n)=0$, then there exists $c_\kappa$ only depending on $k,\overline{n}, c, C$ and $n_{cl}$ such that 
\begin{align*}
&\| (G_{s,c}-[H_s+h\Psi_s])(\cdot,\cdot;x_0,z_0)\|_{L^\infty(\{c\leq |x|\leq C\}\times \{c\leq |z-z_0|\leq C\})}\leq \kappa h^{3/2}
\end{align*}
\end{lemma}

\begin{proof}
The continuous symmetric component, after the change of variable $k^2\tau^2=\lambda-d^2=\lambda-k^2n_h^2+k^2n^2_{cl}= k^2n_{cl}^2-k^2\beta^2$, is written as
\begin{align}\label{eq:Gsymm_form_tau}
\begin{split}
G_{s,c}(x,z;x_0,z_0)=&\frac{1}{2\pi}\int_{0}^\infty v_s(x,k^2\tau^2+d^2) v_s(x_0,k^2\tau^2+d^2)\cdot\\
&\cdot\frac{e^{ik\sqrt{n^2_{cl}-\tau^2}|z-z_0|}}{i\sqrt{n^2_{cl}-\tau^2}}
 \frac{k^2\tau^2}{k^2\tau^2+d^2 \sin^2(h\sqrt{k^2\tau^2+d^2})} d\tau,
\end{split}
\end{align}
where, for $|x|\geq c>h$, we have
\begin{align*}
v_s(x,k^2\tau^2+d^2)&=
\cos(h\sqrt{k^2\tau^2+d^2})\cos(k\tau(|x|-h))\\
&-\sqrt{k^2\tau^2+d^2}\sin(h\sqrt{k^2\tau^2+d^2})\frac{\sin(k\tau(|x|-h))}{k\tau}.
\end{align*}
We will again use the notation $O(h^m)P(\tau)$ to denote any function that, for $h>0$ small can be bounded by $h^m$ times a polynomial in $\tau$, uniformly in the other variables.
Let us start by analyzing the term $h\sqrt{k^2\tau^2+d^2}$. We have
\begin{align*}
h\sqrt{k^2\tau^2+d^2}=\sqrt{h^2d^2+k^2\tau^2}=\sqrt{k^2\overline n^2+h^2k^2(\tau^2-n_{cl}^2)}=k\overline n\sqrt{1+h^2(\tau^2-n_{cl}^2)/\overline n}\\
=k\overline n\left(1+\frac{h^2(\tau^2-n_{cl}^2)}{2\overline n^2}+O(h^4)P(\tau)\right)=k\overline n+\frac{h^2k(\tau^2-n_{cl}^2)}{2\overline n}+O(h^4)P(\tau).
\end{align*}
Since we are assuming that $\sin(k\overline n)=0$, then $\cos(k\overline n)=\pm 1$ and 
\begin{align*}
\sin(h\sqrt{k^2\tau^2+d^2})&=\pm \frac{h^2k(\tau^2-n_{cl}^2)}{2\overline n}+O(h^4)P(\tau),\\
\cos(h\sqrt{k^2\tau^2+d^2})&=\pm 1+O(h^4)P(\tau), \textnormal{ and }\\
\sin^2(h\sqrt{k^2\tau^2+d^2})&=\frac{h^4k^2(\tau^2-n_{cl}^2)^2}{4\overline n^2}+O(h^6)P(\tau).
\end{align*}
We also have that
\begin{align*}
\cos(k\tau(|x|-h))&= \cos(k\tau|x|)\cos(k\tau h)+\sin(k\tau|x|)\sin(k\tau h)\\
&=\cos(k\tau|x|)(1+O(h^2)P(\tau))+k\tau h \sin(k\tau|x|)(1+O(h^2)P(\tau)),\\
&=\cos(k\tau|x|)+k\tau h \sin(k\tau|x|)+O(h^2)P(\tau),\\
h\sin(k\tau(|x|-h))&= h\sin(k\tau|x|)\cos(k\tau h)-h\cos(k\tau|x|)\sin(k\tau h)\\
&= h\sin(k\tau|x|)+O(h^2)P(\tau).
\end{align*}
The estimates above imply that 
\begin{align*}
v_s(x,k^2\tau^2+d^2)&=
\cos(h\sqrt{k^2\tau^2+d^2})\cos(k\tau(|x|-h))\\
&-(h\sqrt{k^2\tau^2+d^2})\frac{\sin(h\sqrt{k^2\tau^2+d^2})}{h^2}\frac{h\sin(k\tau(|x|-h))}{k\tau}\\
&=\pm \left(\cos(k\tau |x|)+ k\tau h \frac{\tau^2+n_{cl}^2}{2 \tau^2}\sin(k\tau |x|) \right)+O(h^2)P(\tau).
\end{align*}
On the other hand
\begin{align*}
\frac{k^2\tau^2}{k^2\tau^2+d^2 \sin^2(h\sqrt{k^2\tau^2+d^2})} 
&=1 -\frac{d^2 \sin^2(h\sqrt{k^2\tau^2+d^2})}{k^2\tau^2+d^2 \sin^2(h\sqrt{k^2\tau^2+d^2})},
\end{align*}
and to estimate the second term in the right hand side we separate in two
intervals. We recall that $h^2d^2=k^2\overline{n}^2+O(h^2)$ while $\sin^2(h\sqrt{k^2\tau^2+d^2})=O(h^4)P(\tau)$ and therefore, for $\tau \in [h^{1/4},\infty)$, 
\begin{align*}
&\frac{d^2 \sin^2(h\sqrt{k^2\tau^2+d^2})}{k^2\tau^2+d^2 \sin^2(h\sqrt{k^2\tau^2+d^2})} =O(h^{3/2})P(\tau),
\end{align*}
For $\tau \in [0,h^{1/4}]$ we use that
\begin{align*}
&\frac{c}{a+c}=\frac{b}{a+b}(1-\frac{a}{a+c}(1-\frac{c}{b})),\quad 0\leq \frac{k^2\tau^2}{k^2\tau^2+d^2 \sin^2(h\sqrt{k^2\tau^2+d^2})}\leq1, \textnormal{ and }\\
&1-\frac{d^2\sin^2(h\sqrt{k^2\tau^2+d^2})}{h^2k^4n_{cl}^4/4}=O(h^{1/2})
\textnormal{ uniformly for } \tau\in[0,h^{1/4}],
\end{align*}
to obtain 
\begin{align*}
\frac{d^2 \sin^2(h\sqrt{k^2\tau^2+d^2})}{k^2\tau^2+d^2 \sin^2(h\sqrt{k^2\tau^2+d^2})}=\frac{h^2k^4n_{cl}^4/4}{k^2\tau^2+h^2k^4n_{cl}^4/4}\Bigg(1-\frac{k^2\tau^2}{k^2\tau^2+d^2 \sin^2(h\sqrt{k^2\tau^2+d^2})}\\
\cdot 
\left(1-\frac{d^2\sin^2(h\sqrt{k^2\tau^2+d^2})}{h^2k^4n_{cl}^4/4}\right)\Bigg)\\
=\frac{h^2k^4n_{cl}^4/4}{k^2\tau^2+h^2k^4n_{cl}^4/4}  (1+O(h^{1/2})) =h\frac{kn_{cl}^2}{2}\left(\frac{h}{(\frac{2\tau}{kn_{cl}^2})^2+h^2}\frac{2}{kn_{cl}^2}\right)(1+O(h^{1/2})),
\end{align*}
uniformly for $\tau\in[0,h^{1/4}]$. In summary, for $\tau\geq 0$
\begin{align*}
&\frac{k^2\tau^2}{k^2\tau^2+d^2 \sin^2(h\sqrt{k^2\tau^2+d^2})} =\\
&1- h\frac{kn_{cl}^2}{2}\left(\frac{h}{(\frac{2\tau}{kn_{cl}^2})^2+h^2}\frac{2}{kn_{cl}^2}\right)\mathbbm{1}_{\{\tau\leq h^{1/4}\}}(1+O(h^{1/2}))+O(h^{3/2})P(\tau).
\end{align*}
In the right hand side, in parenthesis, is the Poisson kernel, for which:
\begin{align*}
&\left|\int_{h^{1/2}}^{h^{1/4}} \!\!\!\!\!f(t) \frac{h}{t^2+h^2}dt\right|\leq O(h^{1/2})\!\!\!\!\!\!\sup_{t\in[h^{1/2},h^{1/4}]}\!\!\!\!\!\!|f(t)|,\quad \int_0^{h^{1/2}} \!\!\!\!\!\!\!f(0) \frac{h}{t^2+h^2}dt = \frac{\pi}{2} f(0) + O(h^{1/2}),\\
& \text{ and }
\left|\int_0^{h^{1/2}} (f(t)-f(0)) \frac{h}{t^2+h^2}dt\right|\leq \frac{\pi}{2}\sup_{t\in[0,h^{1/2}]}|f(t)-f(0)|.
\end{align*}
For
\begin{align*}
f(\tau):= v_s(x,k^2\tau^2+d^2)v_s(x_0,k^2\tau^2+d^2)\frac{e^{ik\sqrt{n_{cl}^2-\tau^2}|z-z_0|}}{i\sqrt{n_{cl}^2-\tau^2}},
\end{align*}
we can check that $|f(\tau)|$ is bounded for $\tau\in[h^{1/2},h^{1/4}]$, that $f(0)=e^{ikn_{cl}|z-z_0|}/in_{cl}+O(h)$ is bounded, and that $\sup_{t\in[0,h^{1/2}]}|f(t)-f(0)|=O(h)$, when $h>0$ is small enough and when $|x|,|x_0|$ and $|z-z_0|$ are in $[c,C]$, with the bounding depending on $k, \overline{n}, c, C$ and $n_{cl}$, but otherwise uniform on the other variables.

With all these estimates on the elements of the integral defining $G_{s,c}$, then equation \eqref{eq:Gsymm_form_tau} can be estimated as
\begin{align*}
&\int_o^\infty f(\tau) \frac{k^2\tau^2}{k^2\tau^2+d^2 \sin^2(h\sqrt{k^2\tau^2+d^2})} d\tau
=\int_0^\infty f(\tau) d\tau -\\
&h \frac{kn_{cl}^2}{2} \int_0^{h^{1/4}} \!\!\!\!f(\tau)\left(\frac{h}{(\frac{2\tau}{kn_{cl}^2})^2+h^2}\frac{2}{kn_{cl}^2}\right)(1+O(h^{1/2}))d\tau+ \int_{h^{1/4}}^\infty f(\tau) O(h^{3/2})P(\tau)d\tau=\\
&\int_0^\infty \cos(k\tau|x|)\cos(k\tau|x_0|)\frac{e^{ik\sqrt{n_{cl}^2-\tau^2}|z-z_0|}}{i\sqrt{n_{cl}^2-\tau^2}}d\tau+\\
&h \int_0^\infty \!\!\!k\tau \frac{\tau^2+n_{cl}^2}{2\tau^2}\sin(k\tau(|x|+|x_0|))\frac{e^{ik\sqrt{n_{cl}^2-\tau^2}|z-z_0|}}{i\sqrt{n_{cl}^2-\tau^2}}d\tau
-h \frac{kn_{cl}^2}{2}\frac{\pi}{2} \frac{e^{ikn_{cl}|z-z_0|}}{in_{cl}}\!+\! O(h^\frac{3}{2}),
\end{align*}
where we used that $\int_0^\infty P(\tau)e^{ik\sqrt{n_{cl}^2-\tau^2}}/\sqrt{n_{cl}^2-\tau^2}d\tau =O(1)$.
Dividing by $2\pi$ and identifying the terms, we complete the proof, since the estimates are uniform for $|x|,|x_0|$ in $[c,C]$ and for $c\leq |z-z_0|\leq C$.

Observe that if we do not want to impose the upper bound $|z-z_0|\leq C$, the proof above allows us to obtain the following estimate for the 0 order term: 
\begin{align*}
&\| (G_{s,c}-H_s)(\cdot,\cdot;x_0,z_0)\|_
{L^\infty(\{c\leq |x|\leq C\}\times \{c\leq |z-z_0|\})}
\leq \kappa h.
\end{align*}
\end{proof}

\subsection{Proof of Theorem \ref{thm:asym_anti}: Asymptotic analysis of the continuous antisymmetric component.}

For the antisymmetric continuous part of the Green's function the analysis is completely analogous.

Recall that we defined $\Phi_a$ and $\Psi_a$ in Theorem \ref{thm:asym_anti}.
The statement (\ref{enum:anti1}) of Theorem \ref{thm:asym_anti} corresponds to the  following lemma.

\begin{lemma}\label{lemma:casym_nr}
Fix $C>c>0$ and assume that $c\leq x_0\leq C$. Also assume that $\cos(k\overline n)\neq0$, then there exists $\kappa $ only depending on $k,\overline{n}, c, C$ and $n_{cl}$ such that 
\begin{align*}
\| (G_{a,c}-[H_a+h\Phi_a])(\cdot,\cdot;x_0,z_0)\|_
{L^\infty(\{c\leq |x|\leq C\}\times \{c\leq |z-z_0|\})}\leq \kappa  h^2.
\end{align*}
\end{lemma}

\begin{proof} The proof is identical to the symmetric case, but observing that the corresponding estimates in this case are
\begin{align*}
hv_a(x,k^2\tau^2+d^2)=& \sgn(x)\Big[ k\overline n\cos(k\overline n)\frac{\sin(k\tau |x|)}{k\tau}\\
&\quad+h\big(\sin(k\overline n)-k\overline n\cos(k\overline n)\big)\cos(k\tau |x|)\Big]+O(h^2)P(\tau),
\end{align*}
and
\begin{align*}
&\frac{2k^3\tau^2}{h^2k^2\tau^2+h^2d^2 \cos^2(h\sqrt{k^2\tau^2+d^2})} =
\frac{2k^3\tau^2}{k^2\overline n^2\cos^2(k\overline n)},
+O(h^2)P(\tau).
\end{align*}
And all the estimates are uniform for $|x|,|x_0|$ in $[c,C]$ and for $|z-z_0|\geq c$.
\end{proof}

The statement (\ref{enum:anti2}) of Theorem \ref{thm:asym_anti} corresponds to the  following lemma.

\begin{lemma}\label{lemma:casym_r}
Fix $C>c>0$ and assume that $c\leq x_0\leq C$. Also assume that $\cos(k\overline n)=0$, then there exists $\kappa $ only depending on $k,\overline{n}, c, C$ and $n_{cl}$ such that 
\begin{align*}
&\| (G_{a,c}-[H_s+h\Psi_a])(\cdot,\cdot;x_0,z_0)\|_
{L^\infty(\{c\leq x\leq C\}\times \{c\leq |z-z_0|\leq C\})}
\leq \kappa h^{3/2},\\
&\text{ and }\\
&\| (G_{a,c}-[-H_s+h\Psi_a])(\cdot,\cdot;x_0,z_0)\|_
{L^\infty(\{-C\leq x\leq -c\}\times \{c\leq |z-z_0|\leq C\})}
\leq \kappa h^{3/2},\\
\end{align*}
\end{lemma}
\begin{proof}
The proof is identical to that of Lemma \ref{lemma:csym_r}, where the only difference is the estimate for the antisymmetric mode
\begin{align*}
v_a(x,k^2\tau^2+d^2)=\pm\sgn(xx_0)\!\left(\cos(k\tau |x|) + k\tau h\frac{3\tau^2-n_{cl}^2}{2\tau^2}\sin(k\tau|x|)\right)\! + O(h^2)P(\tau).
\end{align*}

\end{proof}

\section{Acknowledgments}

E.B. gratefully acknowledges the support from Projets Initiatives de 
Re\-cher\-che Université Grenoble Alpes project 
IRGA2023–IPBC\-–G7H-IRG\-23\-D09. M.C and A.O. gratefully 
acknowledge the support
from ANID-Fon\-de\-cyt 1240200 and ANID Millennium Nucleus ACIP 
NCN19-161. A.O. also acknowledges CMM FB210005 Basal-ANID, 
FONDAP/15110009, ANID-Fon\-de\-cyt 1231404, FONDEF IT23I0095 and 
DO ANID Technology Center DO-\-210001.


\end{document}